\documentclass[12pt]{amsart}
 \usepackage{amsmath,amsthm,amsfonts,amssymb}
\usepackage{verbatim}
\usepackage{latexsym}

\usepackage{amscd}
\usepackage{amssymb}
\usepackage{xcolor}

\newtheorem{theorem}{Theorem}[section]
\newtheorem{lemma}[theorem]{Lemma}
\newtheorem{proposition}[theorem]{Proposition}
\newtheorem{definition}[theorem]{Definition}
\newtheorem{corollary}[theorem]{Corollary}

\newtheorem{remark}[theorem]{Remark}

\newtheorem{thm}{Theorem}[section]

\theoremstyle{definition}

\def\O{\Omega}
\newcommand{\C}{\mathbb C}

\newcommand{\R}{\mathbb R}
\newcommand{\Z}{\mathbb Z}
\newcommand{\N}{\mathbb N}

\def \l{\lambda}
\def \a{\alpha}

\renewcommand{\b}{\beta}
\newcommand{\s}{\mathbf s}

\newcommand{\n}{\mathfrak n}

\def\sp{\text{span}}

\def\k{\mathfrak k}
\def\h{\mathfrak h}

\def \eO{\mathcal O}

\def \rank {\text{\rm rank}}
\def \s{\mathbf s}

\def\Re{\text{\rm Re }}
\def\Im{\text{\rm Im }}
\def \V { \mathcal V}
\def \W { \mathcal W}
\def \Z { \mathcal Z}
\def \E { \mathcal E}
\def \gl {\mathfrak{gl}}

\def\a{\alpha}
\def\b{\beta}

\begin{document}

\title[Classifying wavelet transforms in 3D]{A classification of continuous wavelet transforms in dimension three}

\author{Bradley Currey, Hartmut F\"uhr, Vignon Oussa}
\address{Department of Mathematics and Statistics \\Saint Louis University\\St. Louis. MO 63103\\ U.S.A.}\email{curreybn@slu.edu}

\address{Lehrstuhl A f\"ur Mathematik \\ RWTH Aachen University, D-52056 Aachen\\Germany} \email{fuehr@matha.rwth-aachen.de}

\address{Department of Mathematics \\ Bridgewater State University\\ Bridgewater, MA 02324 \\U.S.A.}\email{vignon.oussa@bridgew.edu}

\begin{abstract}
This paper presents a full catalogue, up to conjugacy and subgroups of finite index, of all matrix groups $H < {\rm GL}(3,\mathbb{R})$ that give rise to a continuous wavelet transform with associated irreducible quasi-regular representation. For each group in this class, coorbit theory allows to consistently define spaces of sparse signals, and to construct atomic decompositions converging simultaneously in a whole range of these spaces.  
As an application of the classification, we investigate the existence of compactly supported admissible vectors and atoms for the groups. 
\end{abstract}

 \maketitle
%


\section{Introduction}

Continuous wavelet transforms in higher dimensions have been studied since \cite{Mu}. The construction departs from the choice of a closed matrix group $H< {\rm GL}(d,\mathbb{R})$. The affine group generated by $H$ and the group of translations is the semi-direct product group $G = \mathbb{R}^d \times H$. This group acts unitarily on functions $f \in {\rm L}^2(\mathbb{R})$ by the quasi-regular representation
\[
 \pi(x,h) f(y) = |{\rm det}(h)|^{-1/2} f(h^{-1}(y-x))~.
\] Given any pair of vectors $\psi,f \in {\rm L}^2(\mathbb{R}^d)$, the {\bf continuous wavelet transform} of $f$ with respect to $\psi$ is the function 
\[
 \mathcal{W}_\psi f : G \to \mathbb{C}~,~ (x,h) \mapsto \langle f, \pi(x,h) \psi \rangle~.
\] The continous wavelet transform associated to $\psi$ is the linear operator $\mathcal{W}_\psi: f \mapsto \mathcal{W}_\psi f$, and $\psi$ is called a {\bf wavelet} or {\bf analyzing vector}. We call $\psi$ {\bf admissible} if $\mathcal{W}_\psi$ is a nonzero scalar multiple of an isometry into ${\rm L}^2(G)$, i.e., if the equality
\[
 \| f \|_2^2 = \frac{1}{c_\psi} \int_H \int_{\mathbb{R}^d} \left| \mathcal{W}_\psi f (x,h) \right|^2 dx \frac{dh}{|{\rm det}(h)|} 
\] holds for all $f \in {\rm L}^2(\mathbb{R}^d)$. This isometry property allows to invert the wavelet transform by its adjoint, which leads to the weak-sense inversion formula
\[
 f =  \frac{1}{c_\psi} \int_H \int_{\mathbb{R}^d}  \mathcal{W}_\psi f (x,h) ~\pi(x,h) \psi~ dx \frac{dh}{|{\rm det}(h)|} 
\]

We call $H$ {\bf  admissible} if there exist admissible vectors $\psi \in {\rm L}^2(\mathbb{R}^d)$, and {\bf irreducibly admissible} if the quasi-regular representation is in addition irreducible. 
Admissible groups have been studied extensively, see \cite{BeTa,Fu96,LaWeWeWi,GrKaTa,FuMa,Fu_LN,Fu_Calderon,DaKuStTe,Bruna_etal,ACDO} for a small sample of the existing literature. 
Nonabelian versions, where $\mathbb{R}^d$ is replaced by a suitably chosen noncommutative Lie group, were also considered \cite{
CuMcN,Currey07,GeMa,Ou13}.
A particular class of irreducibly admissible matrix groups gained particular attention, the so-called {\bf shearlet group} acting in dimension two \cite{DaKuStTe} and its higher-dimensional relatives \cite{DaStTe12,DaHaTe12,CzaKi12,Fu_RT}.

The interest of irreducibly admissible dilation groups is strengthened by the realization  that coorbit space theory, as developed in \cite{FeiGr0,FeiGr1,FeiGr2}, applies to these groups and their associated continuous wavelet transform. This observation had been made early on for the similitude groups in \cite{FeiGr1}, was later shown for the shearlet cases in \cite{DaKuStTe,DaStTe12,DaHaTe12}, and finally established in full generality in \cite{Fu_coorbit}. Since then, the theory of wavelet coorbit spaces has been further developed in \cite{Fu_atom,Fu_RT,FuVo,VoDiss}; see also \cite{FeVo} for an overview of the recent developments and their historical roots. 

Coorbit space theory provides a unified approach to wavelet approximation theory, that allows to consistently define function spaces in terms of wavelet coefficient decay. Here consistency refers to the fact that for functions in these spaces, the fast wavelet coefficient decay is guaranteed for {\em all} choices of analyzing vectors taken from a well-understood class of nice wavelets. This property emphasizes that the approximation-theoretic properties of group-theoretically defined wavelet systems ultimately depend on the way the group $G$ generates the system from a single mother wavelet, rather than on the specific choice of the mother wavelet.  In other words, it puts the focus of the analysis of wavelet system on the role of the dilation group $H$. It is the main purpose of this paper to provide an overview of the available choices of $H$ in dimension three.

\section{Aims and main result of this paper}

In this paper we set out to determine all irreducibly admissible matrix groups $H< {\rm GL}(3,\mathbb{R})$ up to conjugacy. In effect, we will only classify the connected components $H_0< H$ of the identity in $H$. However, from this subgroup to the full group is only a small step: By \cite[Remark 4]{Fu98}, $H_0$ is an open subgroup of finite index. 

There are several motivations for studying such a classification. The chief interest of obtaining a full list of candidate groups is to get an overview of the variety of available group actions.
The three-dimensional case is arguably of interest to applications, due to the increasing use made of voxel data provided by imaging applications such as computerized tomography or MRI. Judging by the existing digital implementations of wavelet-like transforms such as curvelets and shearlets in the three-dimensional setting \cite{CaDeDoYi,NeLa}, the task of numerically implementing a given continuous wavelet transform in dimension three seems challenging, but not insurmountable. Hence we believe that the list of available options presented below is potentially useful in applications.  

In such a description of all possible cases, it is however advisable to use suitable equivalence relations as a means of reducing the complexity of the general picture. In our context, conjugacy is the natural choice for such an equivalence relation: If two matrix groups $H_1,H_2$ are conjugate, say $H_1 = g H_2 g^{-1}$ for some invertible matrix $g$, then it is very straightforward to relate the associated wavelet systems via a dilation by $g$. Indeed, if $\psi_2$ is an admissible vector for $H_2$, and $f_2 \in {\rm L}^2(G)$ is arbitrary, then defining 
\[
 \psi_1(x) = |{\rm det}(g)|^{-1/2} \psi_2(g^{-1}(x)~,f_1(x) = |{\rm det}(g)|^{-1/2} f_2(g^{-1}(x)
\] 
we obtain for $h_i \in H_i (i=1,2)$ related via $h_1 = gh_2 g^{-1}$ and $x \in \mathbb{R}^d$ that 
\[
 \mathcal{W}_{\psi_1}^{H_1} f_1(x,h_1) = \mathcal{W}_{\psi_2}^{H_2}(g^{-1}x, h_2)~.
\]
Here we temporarily used $H_1,H_2$ as superscripts to differentiate the wavelet transforms associated to $H_1,H_2$, respectively. Hence one can consider the wavelet transforms over $H_1$ and $H_2$ as instances of the same transform, viewed with respect to different linear coordinate systems on $\mathbb{R}^d$. In particular, any property of interest pertaining to wavelets arising from the action of $H_1$ will translate to an analogous property of wavelets arising from $H_2$, and it is clearly sufficient to study just one of the groups. More generally speaking, for most questions of interest it will be sufficient to study just one representative from each full conjugacy class. 

 Our chief example of a property that can be checked by focusing on the representatives modulo conjugacy is the question whether there exist compactly supported admissible vectors and atoms for a given wavelet transform. It is currently not known whether every irreducibly admissible matrix group has compactly supported admissible vectors, and for the much smaller class of atoms (also called ``better vectors'' in the coorbit scheme \cite{FeiGr1}), this existence question is even harder to answer. We prove, via a case-by-case treatment of all possible candidates, that every irreducibly admissible matrix group in dimension three has compactly supported admissible vectors, and that the overwhelming majority of cases even allows the existence of compactly supported atoms.  

In the context of generalized wavelet analysis, some classification results are already known. For instance, the analogous search for all irreducibly admissible $H< {\rm GL}(2,\mathbb{R})$ groups in dimension two (or at least, of their connected components) produces the following list, complete up to conjugacy \cite{FuCuba}:
\begin{itemize}
 \item The {\bf diagonal group}, given by 
 \[
  H  = \left\{ \left[ \begin{array}{cc} a & 0 \\ 0 & b \end{array} \right]~:~a>0, b>0 \right \}~;
 \]
 \item the {\bf similitude group} 
  \[
  H  = \left\{ \left[ \begin{array}{cc} a & b \\ -b & a \end{array} \right]~:~a^2 + b^2 \not= 0 \right \}~;
 \]
 \item and the {\bf shearlet groups} defined as 
 \[
  H  = \left\{ \left[ \begin{array}{cc} a & b \\ 0 & a^c \end{array} \right]~:~a>0, b \in \mathbb{R} \right \}~.
 \] Here $c \in \mathbb{R}$ can be chosen arbitrarily. 
\end{itemize}
Existence of atoms for each of these cases (and thus for all irreducibly admissible matrix groups acting in dimension two) was shown in \cite{Fu_atom}. The analytic properties of the associated wavelet transforms were explored in the thesis \cite{Miro}.
Further (partial) classification results exist for {\em abelian} irreducibly admissible groups in arbitrary dimensions \cite{Fu98}, and for the class of generalized shearlet dilation groups \cite{Fu_RT} . In a similar spirit to our paper, the papers \cite{AlBaDeMDeV,AlBaDeMDeVMa} classified a class of subgroups of the metaplectic group that allow an inversion formula under the metaplectic representation. We are not aware of further sources dealing with the systematical construction and classification of irreducibly admissible matrix groups or certain subclasses thereof. 

The following theorem summarizes the results of our paper. Throughout the paper, we use $\mathbb{R}^* = \mathbb{R} \setminus \{ 0\}$. The terminology appearing in part (c) will be introduced in Section \ref{sect:atom}.
\begin{theorem}
 \label{thm:main}
 Let $H< {\rm GL}(3,\mathbb{R})$ be irreducibly admissible.  
\begin{enumerate}
 \item[(a)] Its connected component is conjugate to precisely one group from the following list: 
\begin{enumerate}
 \item[(1)] {\bf The abelian groups.}
 \begin{enumerate}
 \item[(1.a)] 
\[ H_0 = \left\{ \left[ \begin{array}{ccc} a_1 &  & 0 \\ 0 & a_2 & 0 \\ 0 & 0 & a_3 \end{array} \right] : a_i >0 
\right\}  .\]

\item[(1.b)] 
\[ H = \left\{ \left[ \begin{array}{ccc} a & 0 & 0 \\ 0 & b \cos t & -b\sin t \\ 0 & b\sin t & b\cos t \end{array} \right] : a> 0,b> 0, t \in \mathbb{R}
 \right\}  .\]

\item[(1.c)] 
\[ H_0 = \left\{ \left[ \begin{array}{ccc} a & b & 0 \\ 0 & a & 0 \\ 0 & 0 & c \end{array} \right] : a>  0,c>  0, b \in \mathbb{R}
\right\}  .\]

\item[(1.d)] 
\[ H_0 = \left\{ \left[ \begin{array}{ccc} a & b & c \\ 0 & a & b \\ 0 & 0 & a \end{array} \right] : a>  0 , b, c \in \mathbb{R}
\right\}  .\]

\item[(1.e)] 
\[ H_0 = \left\{ \left[\begin{array}{ccc} a & b & c \\ 0 & a & 0 \\ 0 & 0 & a \end{array} \right] : a>0, b, c \in \mathbb{R}
\right\}  .\]
\end{enumerate}
 \item[(2)] {\bf The nonabelian, solvable groups.}
 \begin{enumerate}
 \item[(2.a)]
$$
H_0 = \left\{ \left[ \begin{matrix} e^{as}a &t_1&t_2\\0&e^{as}\cos(bs)&-e^{as}\sin(bs)\\0&e^{as} \sin(bs) & e^{as} \cos(bs) \end{matrix} \right] : s, t_1, t_2 \in \R\right\},
$$
with $(a,b) \in \mathbb{R} \times \mathbb{R}^\ast$. 
\item[(2.b)]  
$$
H_0 = \left\{ \left[ \begin{matrix} e^s&t_1&t_2\\0&\cos(r)&-\sin(r)\\0& \sin(r) & \cos(r) \end{matrix} \right] : r, s, t_1, t_2 \in \R\right\}.
$$
\item[(2.c)] 
$$
H_0 = \left\{ \left[ \begin{matrix} e^{as}&t_1&t_2\\0&e^{s}\cos(r)&-e^{s}\sin(r)\\0&e^{s} \sin(r) & e^{s} \cos(r) \end{matrix} \right] : r, s, t_1, t_2 \in \R\right\},
$$
with $a \in \mathbb{R}^*$.
\item[(2.d)] 
$$
H_0 = \left\{ \left[\begin{matrix} a^\l & t_1&t_2\\ 0& a^{\l}&0\\ 0&0&a^{\l-1}\end{matrix}\right] : a > 0, t_1, t_2 \in \R \right\},
$$
with $\lambda \in \mathbb{R}^*$.
\item[(2.e)]
$$
H_0 = \left\{ \left[\begin{matrix} a^\l & t_1 &t_2\\ 0& a^{\l-c} &0\\ 0&0&a^{\l-1}\end{matrix}\right] : a > 0, t_1, t_2 \in \R\right\},
$$
with $(\lambda,c) \in \mathbb{R} \times ([-1,1] \setminus \{ 0 \})$.
\item[(2.f)]
$$
H_0 = \exp \h = \left\{ \left[\begin{matrix} a^\l& t_1&t_2  \\ 0& a^\l &a^\l \ln a \\ 0&0&a^\l \end{matrix}\right] : a >0, t_1, t_2 \in \R\right\},
$$ with $\lambda \in \mathbb{R}^*$.
\item[(2.g)]
$$
H_0 = \exp \h = \left\{ \left[\begin{matrix} a^\l& t_1&t_2  \\ 0& a^\l &t_1 + a^\l \ln a \\ 0&0&a^\l \end{matrix}\right] : a >0, t_1, t_2 \in \R\right\},
$$
with $\lambda \in \mathbb{R}^*$. 
\item[(2.h)]
$$
H_0 = \left\{ \left[\begin{matrix} a^\l& t_1 &t_2\\ 0& a^{\l-1}&\delta  a^{\l-1}\ln a\\ 0&0&a^{\l-1}\end{matrix}\right] : a>0,  t_1, t_2 \in \R\right\}.
$$ with $(\lambda, \delta) \in \mathbb{R}^\ast \times \mathbb{R}$. 
\item[(2.i)]
$$
H_0 = \left\{ \left[\begin{matrix} a^\l& t_1 &t_2\\ 0& a^{\l-1}& t_1 a^{\l-2}\\ 0&0&a^{\l-2}\end{matrix}\right] : a > 0, t_1, t_2 \in \R\right\},
$$
with $\lambda \in \mathbb{R}^*$. 
\item[(2.j)]
$$
H_0 = \left\{ \left[\begin{matrix} a^{\l+1} b^c & t&0\\ 0& a^\l b^c&0\\ 0&0&b \end{matrix}\right] : a>0, b>0, t \in \R\right\}
$$
with $(\lambda,c) \in \mathbb{R}^* \times (\mathbb{R} \setminus \{ 1 \})$.
\item[(2.k)]
$$
H_0 = \left\{ \left[\begin{matrix} a b & t&0\\ 0& b&0\\ 0&0&a^\l  \end{matrix}\right] : a>0, b>0, t \in \R\right\}
$$ with $\lambda \in \mathbb{R}^*$.
\item[(2.l)]
$$
H_0 = \left\{\left[\begin{matrix} b& t &-\nu_1\ln a + \nu_2b \ln b \\ 0& ab&0\\ 0&0&b\end{matrix}\right] : a>0, b>0, t \in \R\right\}.
$$ with $(\nu_1,\nu_2) \in \mathbb{R}^* \times \mathbb{R}$. 
\item[(2.m)]
$$
H_0 = \left\{ \left[\begin{matrix} ab & t &0\\ 0& b&0\\ 0&0&b^\l\end{matrix}\right] : a>0, b>0, t \in \R\right\}.
$$ with $\lambda \in \mathbb{R}^*$.
\item[(2.n)] 
$$
H_0 = \left\{ \left[\begin{matrix} ab&0 &t\\ 0& b &\delta_1  \ln a  + \delta_2 b \ln b\\ 0&0&b\end{matrix}\right] : a>0, b>0, t \in \R\right\}.
$$ with $(\delta_1,\delta_2) \in \mathbb{R}^2 \setminus {(0,0)}$.
 \end{enumerate}
\item[(3)] {\bf The nonsolvable groups.}
\begin{enumerate}
\item[(3.a)] $H = \mathbb{R}^+ \times SO(3)$. 
\end{enumerate}
\end{enumerate}
Every group $H_0$ appearing in this list is contained in an irreducibly admissible matrix group $H$.
\item[(b)] $H$ has compactly supported admissible vectors. 
\item[(c)] If the connected component of $H$ is not conjugate to one of the groups from case (2.l), then $H$ has compactly supported atoms $\psi \in \mathcal{B}_{v_0}$, for any polynomially bounded weight $v_0$ on $G = \mathbb{R}^3 \rtimes H$. 
\end{enumerate}
\end{theorem}

The classification of the solvable cases is achieved in Section \ref{sect:solv}, whereas the classification of the non-solvable cases can be found in Section \ref{sect:nonsolv}. Parts (b) and (c) are proved in Section \ref{sect:atom}.

We count 7 isolated cases (among them the 5 abelian ones), 7 one-parameter and 6 two-parameter families, with the overwhelming majority of examples provided by the solvable cases. The cases (1.e), (2.d) and (2.e) correspond to the standard shearlet groups \cite{DaKuStTe}, whereas the case (2.i) contains the Toeplitz shearlet group \cite{DaHaTe12}. We note that some of the parametrized families given above can be merged, e.g. (1.e), (2.d) and (2.e), which would lead to a somewhat smaller list of cases. However, the current list reflects the construction method for the different cases and thus enables an easier tracing of the argument proving the classification, which is why we have refrained from streamlining the presentation in the theorem.

\section{Preliminaries} 

\subsection{Properties of irreducibly admissible matrix groups}

Our classification of irreducibly admissible matrix groups $H$ departs from a {\em characterization} of said groups, that uses the dual action of $H$. This group action is given by the mapping
$H \times \mathbb{R}^d \ni (h, \xi) \mapsto h^{-T} \xi$, where $h^{-T}$ denotes the inverse of the transpose of $h$. Now irreducibly admissible matrix groups are characterized as follows \cite{Fu_Calderon},  \cite[Corollary 5.15]{Fu_LN}.
\begin{theorem}
Let $H< {\rm GL}(d,\mathbb{R})$ denote a closed matrix group.
\begin{enumerate}
\item[(a)]  $H$ is irreducibly admissible if and only if there exists an open, conull orbit $\mathcal{O} \subset \mathbb{R}^d$ with the additional property that for all $\xi \in \mathcal{O}$, the associated stabilizer
 \[
  H_{\xi} = \{ h \in H: h^{-T} \xi = \xi \}
 \] is compact.
 \item[(b)] If $H$ is admissible, then $G = \mathbb{R}^d \rtimes H$ is nonunimodular.
 \end{enumerate}
\end{theorem}
Note that the condition in part (b) is equivalent to the statement that the functions $\Delta_H$ and $|{\rm det}|$ do not coincide; here $\Delta_H$ is the modular function of $H$. 

By \cite[Remark 4]{Fu98}, property (a) implies that the connected component $H_0< H$ acts with open orbits and (almost everywhere) compact stabilizers. Furthermore, the fact that the stabilizers in $H$ associated to the open orbits are compact implies that $H/H_0$ must be finite. Also by \cite[Remark 4]{Fu98}, the dimension of admissible matrix groups $H < {\rm GL}(d,\mathbb{R})$ lies between $d$ and $d(d+1)/2$, which for $d=3$ means that ${\rm dim}(H) \in \{ 3, 4 \}$. To summarize, our search in the following concentrates on connected matrix groups $H_0< {\rm GL}(3,\mathbb{R})$ with open orbits and associated compact stabilizers; necessarily of dimension $3$ or $4$. 
Our proof strategy consists in carefully studying the dual action of connected groups $H< {\rm GL}(3,\mathbb{R})$, and to extract conjugation-invariant features of this action that allow to systematically separate the general class into a various subcases giving rise to the classification in Theorem \ref{thm:main}. For this purpose, we rely on various notions from differential geometry and Lie group theory. These are introduced next. 

\subsection{Linear actions of Lie groups}

Let $\V$ be a finite-dimensional vector space over $\R$, and let $H$ be a locally compact topological group. A {\it representation} of $H$ in $\V$ is a continuous homomorphism $\a: H\rightarrow GL(\V)$; a representation of $H$ can be regarded as a linear action of $H$ on $\V$. If $\W$ is the complexification of $\V$, then $\a$ defines a linear action of $H$ on $\W$ in the obvious way. Given a representation $\a$ of $H$ in $V$, the Cartesian product $\V \times H$ is endowed with the semi-direct product operation defined by
$$
(x,h) \cdot (x',h') = (x+ \a(h)x', hh')
$$
and the resulting group is denoted by $\V \rtimes_\a H$. Let $\a : H \rightarrow GL(\V)$ and $\b : K \rightarrow GL(\W)$ be representations of groups $H$ and $K$. We say that $\a$ and $\b$ are isomorphic if there is a group isomorphism $c : H \rightarrow K$ and a vector space isomorphism $g : \V \rightarrow \W$, such that 
$$
\b(c(h))\circ g = g \circ \a(h).
$$
Note that in this situation $(v,h) \mapsto (gv, c(h))$ is an isomorphism of the corresponding semi-direct products. 
Given a representation $\a$, its dual representation $\a^* : H \rightarrow GL(\V^*)$ is defined by $\a^*(h) f = f\circ \a(h)^{-1}$. 
If $\a$ and $\b$ are isomorphic representations, then $\a^*$ and $\b^*$ are isomorphic.\\

The groups of interest in this paper are matrix groups. We identify an element $X\in \gl(d,\R)$ with a $d\times d$ real matrix in the usual canonical way. With the canonical inclusion $\R^d \subset \C^d$, we have $\gl(d,\R) \subset \gl(d,\C)$. The conjugacy class of an element $X \in \gl(d,\R)$ is the set $\{ s X s^{-1} : s \in GL(d,\C)\}$; note that  $X'$ belongs to the conjugacy class of $X$ if and only if $X'$ is the matrix representation of $X$ with respect to some basis of $\C^n$. 

Let $H$ be a subgroup of $GL(d,\R)$. If $\V = \R^d$, then $\a =  \mathbf{id}_H$ is called the standard representation of $H$. The representation of interest in this paper is the dual of the standard representation: identifying $\R^d$ with its dual via the standard inner product, the dual of the standard representation is given by $h \mapsto (h^T)^{-1}$, and is called the transpose representation. If $H$ and $K$ are subgroups of $GL(d,\R)$, then the standard representations of $H$ and $K$, or their duals, are isomorphic if and only if $H$ and $K$ are conjugate.

If $H$ is a {\it closed} subgroup of $GL(d,\R)$, then $H$ is a Lie subgroup of $GL(d,\R)$. For the present purposes, this means that each of the following hold. 

\medskip
\noindent
(1) There is a subspace $\h$ of the space $\mathfrak{gl}(d,\R)$ of all endomorphisms of $\R^d$ that is closed under the usual bracket operation $[X,Y] = XY - YX$, and such that 
$$
\h = \{ X \in \mathfrak{gl}(d,\R) : \exp \R X \subset H\}
$$
where 
$$
\exp X = I + X + \frac{1}{2!} X^2 + \frac{1}{3!} X^3 + \cdots;
$$
$\h$ is called the Lie algebra of $H$. 

\medskip
\noindent
(2) If $n = \dim \h$, then there is an open neighborhood $U$ of the identity $1$ in $H$, and an open neighborhood $V$ of $0 \in \h$, such that $\exp |_V : V\rightarrow U$ is a homeomorphism. 
More generally, if $K$ is a closed subgroup of $H$ with Lie algebra $\k$, then there is an open neighborhood $U$ of $K$ in the quotient space $H/K$, and an open neighborhood $V$ of a complementary space $\mathfrak{s}$ for $\k$ in $\h$, such that the map $X \mapsto \exp X K$ is a homeomorphism of $V$ onto $U$.

\medskip
\noindent
(3) If $\a : H \rightarrow GL(\V)$ is a representation of $H$, then there is a unique Lie algebra homomorphism $d\a : \h \rightarrow \mathfrak{gl}(\V)$ such that $\a(\exp X) = \exp d\a(X)$ holds for all $X \in \h$. One has 
$$
d\a(X) = \frac{d}{dt} \a(\exp tX) \Big|_{t = 0}, \ \ X \in \h
$$
and the Lie algebra of $\ker \a$ is the kernel of $d\a$ in $\h$.

\medskip
\noindent
A Lie algebra homomorphism of $\h$ into $ \mathfrak{gl}(\V)$ is called a Lie algebra representation of $\h$. Just as a representation of $H$ is a continuous linear group action of $H$, a Lie algebra homomorphism of $\h$ into $ \mathfrak{gl}(\V)$ is said to define a linear action of the Lie algebra $\h$ on $\V$.  

\

If $H$ and $K$ are closed subgroups of $GL(d,\R)$ and $c : H \rightarrow K$ is a continuous isomorphism, then $dc$ is an isomorphism of their Lie algebras, and if $\l$ is a linear form on $\h$ then $\l^c := \l(dc^{-1} \cdot)$ is a linear form on $\k$ and $\l \mapsto \l^c$ is a vector space isomorphism of $\h^*$ and $\k^*$. In particular if $K = sHs^{-1}$ for some $s \in GL(n,\R)$, then the Lie algebras $\h$ and $\k$ of $H$ and $K$ satisfy $\k = s\h s^{-1}$. 

\

From now on $H$ will denote a closed subgroup of $GL(d,\R)$ with Lie algebra $\h$, acting on $\R^d$ (and $\C^d$) by the transpose representation. A subspace $\Z$ of $\R^d$ is $H$-invariant if $h^T\Z  = \Z$ for all $h \in H$, and an invariant subspace $\Z$ is called $H$-irreducible if $\Z \ne \{0\}$ and the only proper subspace of $\Z$ that is invariant is $\{0\}$. Similar notions of invariance and irreducibility exist for representations of a Lie algebra. 


\

Denote by $H_0$ the connected component of the identity in $H$; then $H_0$ is generated by elements of the form $\exp X, X \in \h$. It follows that  $\Z$ is $H_0$-invariant if and only if $\Z$ is $\h$-invariant, and $\Z$ is $H_0$-irreducible if and only if $\Z$ is $\h$-irreducible. 

\

Given an invariant subspace $\Z$ of $\R^d$, let $\pi :\V \rightarrow \R^d / \Z$ be the canonical map; $H$ acts linearly on the quotient $\R^d / \Z$ by 
$$
h \cdot \pi(x) = \pi((h^T)^{-1}x), \ \ h \in H, x \in \C^d;
$$
there is an associated quotient action of $\h$ as well. The complexification of the quotient space $\R^d / \Z$ is $\C^d / \Z_\C$. 
Given any linear action of $\h$ in a real vector space $\V$, 
a complex linear form $\l$ on $\h$ is a {\it characteristic function} for this action if the subspace 
$$
\E_\l = \{ u \in \mathcal V_\C  : X \cdot u = \l(X) u, \forall X \in \h \}
$$
is non-zero. The non-zero elements of $\E_\l$ are called characteristic vectors for $\l$. Suppose that $\l$ is a characteristic function for a linear action, and choose a characteristic vector $e$ for $\l$. If $\l$ is real-valued, then we can choose $e \in \V$ and $\Z = \R e$ is $\h$-invariant. If $\l$ is not real-valued, then write $e = e_1 + i e_2$ with $e_1, e_2 \in V$, and put $\Z = \sp \{ e_1, e_2\}$. Then $\dim_\R \Z = 2$, $\Z$ is $\h$-invariant, and the linear action of $\h$ in $\Z$ has the characteristic functions $\l$ and $\overline{\l}$. 
Moreover, the subspace $\Z$ defined in each case above is both invariant and irreducible. 

A complex linear form on $\h$ is a {\it weight} if there is an invariant subspace $\Z$ of $\V_\C$ such that $\l$ is a characteristic function for the quotient action of $\h$ on $\V_\C / \Z$.

\

 The following is easily verified. 
 
 \begin{lemma} \label{weight kernel}  Let $\h$ act linearly on a real vector space $\V$. 
 
 \begin{itemize}
 
 \item[{\rm (a)}]  For any weight $\l$, $[\h,\h] \subset \ker \l$. 
 
 \item[{\rm (b)}] If $\l$ is a weight, then $\overline{\l}$ is also a weight, and $\E_{\overline{\l}} = \overline{\E}_\l$.  
 
 \item[{\rm (c)}] If there are $d$ independent weights, then $\h$ is commutative.

 \end{itemize} 
 
 \end{lemma}

The linear action of primary interest here is the transpose action, that is, the dual of the standard representation. With $\R^d$ identified with its dual via the dot product, we have 
$$
h \cdot v = (h^T)^{-1} v, \ \ \ h \in H, v \in \R^d
$$
and the corresponding action of the Lie algebra $\h$ is 
$$
X \cdot v = - X^T v, \ \ \ X \in \h, v \in \R^d.
$$
For $v \in \R^d$, let $H(v) $ be the stabilizer of $v$ for the transpose action, and let $\h(v ) = \{ X \in \h : X^T v = 0\}$) be the annihilator of $v$ in $\h$. Then $H(v)$ is a closed subgroup of $H$ (though not necessarily connected), and $\h(v)$ is a Lie subalgebra of $\h$. There is a natural continuous bijection from the quotient space $H/H(v)$ onto the $H^T$-orbit of $v$ given by $hH(v) \mapsto (h^T)^{-1} v$, and this map is a homeomorphism if and only if the orbit is a locally compact subset of $\R^d$. We leave the proof of the following to the reader. 

\begin{lemma}  \label{weights standard transpose} A linear form $\l$ on $\h$ is a weight for the transpose action of $\h$ if and only if $- \l$ is a weight  for the standard action.

\end{lemma}

Define $\phi = \phi_v : H \rightarrow \V$ by $\phi(h) = (h^T)^{-1} v$. The differential of $\phi$ at the identity $e \in H$ is a map $d\phi_e  : \h \rightarrow \V$ given by 
$$
 \frac{d}{dt} \Big|_{t = 0}f(\exp tX^T)^{-1}v) = \langle \bigtriangledown f(v) , d\phi_e(X) \rangle
$$
where $f$ is any smooth function on $\V$. The following is an immediate consequence of the chain rule.

\begin{lemma} \label{ann} One has $d\phi_e(X) = -X^Tv$. Hence 
$
\h(v ) = \ker  d\phi_e
$
and coincides with the Lie algebra of $H(v)$.
\end{lemma}

Let $\eO $ be an $H^T$-orbit in $\R^d$. Given $v \in \eO$ and $w = hv$, then $\h(w) = h \h(v) h^{-1}$, and hence the map $v \mapsto \dim ( \h/\h(v))$ is constant on $\eO$. In fact it can shown that $\eO$ is a submanifold of $\R^d$ whose tangent space at each $v \in \eO$  is the image of $d \phi_v$, so that $\dim \eO = \dim \h / \h(v)$. In particular, for any $v \in \R^d$, $\dim \h/ \h(v) \le d$. 

The following lemma translates the open orbit property for $H$ to a condition on $\mathfrak{h}$.




\begin{lemma} \label{open orbit} Let $\eO = H^Tv \subset \R^d$. Then $\dim \h / \h(v) = d $ if and only if $\eO$ is open. 
\end{lemma}

\begin{proof} Suppose that $\dim \h / \h(v) = d$. Then for each $v \in \eO$,  $(d\phi_v)_e$ is an isomorphism of $\h / \h(v) $ onto $\R^d$. Now the inverse mapping theorem says that there is an open neighborhood  $U$ of $e$ in $H$ and an open neighborhood $E$ of $v$ in $\R^d$ such that $\phi|_U$ is a homeomorphism of $U$ onto $E$. In particular the orbit $\eO$ contains an open neighborhood of $v$. 

Conversely, suppose that $\eO$ is open and fix $v \in \eO$. We first show that the quotient space $H/H(v)$ is homeomorphic with $\eO$. Let $K$ be a compact neighborhood of $H(v)$ in $H/H(v)$. The canonical image of $K$ in $\eO$ is a compact continuous image of $K$. Now $H/H(v)$ is a countable union of translates of $K$, and hence $\eO$ is a countable union of translates of the image of $K$. It follows that the image of $K$ has non-empty interior in $\eO$, hence in $\R^d$. We conclude that there is an open neighborhood $U$ of $H(v)$ in $H / H(v)$ and an open neighborhood $V$ of $v \in \eO$ that are homeomorphic. But now any open subset of $H/H(v)$ is a union of translates of some open subset $U_0$ of $U$, and hence its image in $\eO$ is the union of the same translates of the image $V_0$ of $U_0$. Since $V_0$ is open this shows that the canonical map $H /H(v)$ is open, hence a homeomorphism. 

By the Lie group property (2) above, we have a subspace $\mathfrak{e} $ of $\h$, having dimension $\dim(\h/\h(v))$, that contains an open subset homeomorphic with an open subset of $ \R^d$, so $\dim(\h/\h(v)) = d$.

\end{proof}

Put 
$$
\O = \{ v \in \V : \dim \h / \h(v) = d\}.
$$
Fixing a basis $\{X_j : 1 \le j \le n\} $ for $\h$ and $\{ e_i : 1 \le i \le d\}$ for $\R^d$, then $\O = \{ v \in \V : \rank M(v) = d\}$ where $M(v)$ is the $d\times n$ matrix 
$$
M(v) = \left[  \langle e_i, X_j^Tv\rangle\right].
$$
Thus $\O$ is Zariski-open in $\R^d$. Of course there are many examples of Lie subalgebras $\h$ of $\gl(d,\R)$ for which $\O$ is empty. If $\O$ is nonempty, then we say that $\h$ is an open orbit subalgebra of $\gl(d,\R)$. If $H$ is a closed subgroup of $GL(d,\R)$ with Lie algebra $\h$, then the transpose representation of $H$ produces open orbits if and only if $\h$ is an open orbit subalgebra. In this case we say that $H$ is an open orbit group. 

Lemmas \ref{weights standard transpose}, \ref{ann}, and \ref{open orbit} hold for any representation of a Lie algebra in a real vector space. Given a real subspace $\V$ that is invariant for the transpose action, we will also consider the restriction of the transpose action to $\V$, as well as the quotient action on $\R^d / \V$.




\

\section{Classifying the solvable groups}

\label{sect:solv}

The classification of the solvable connected subgroups $H< {\rm GL}(3,\mathbb{R})$ acting with open orbits and associated compact stabilizers constitutes a large portion of this paper, and they contribute all cases except one in the list from Theorem \ref{thm:main}. Let us briefly outline the arguments and results of this paper, and how they combine to provide the major part of the proof of \ref{thm:main}(a). We proceed by analyzing weights on the Lie algebra that are associated to the dual action. While one such weight, associated to the dual action of a particular matrix group $H$, is not a conjugation invariant feature of said group, the {\em existence} of weights with certain additional properties turns out to be such a feature. Similarly, the subspace of nilpotent matrices contained in the Lie algebra of $H$ is not a conjugation invariant, but its dimension is. These features of $H$, or rather its Lie algebra, allows to separate the class of solvable groups into several subclasses, each of which can 
then be 
classified further, to the point of computing a list of representatives modulo conjugacy. More precisely, we obtain the following subclasses:
\begin{itemize}
 \item $H$ is abelian; see Proposition \ref{prop:ab}. This contributes the cases (1.a)--(1.e). Admissible groups containing these groups are given in \ref{prop:ab_nc}.
 \item $H$ is nonabelian, and has a nonreal weight; see Proposition \ref{prop:nrw}.  This contributes the cases (2.a)--(2c), with the admissible groups containing the connected ones provided by Corollary \ref{cor:nrw_nc}.
 \item $H$ is nonabelian, all its weights are real-valued, and the dimension of the nilpotent part is two; see Proposition \ref{prop:rw_np2}. This contributes the cases (2.d)--(2.i); note that case (3) from the Proposition corresponds to the cases (2.f) and (2.g) from Theorem \ref{thm:main}. The admissible groups containing the connected ones are provided by Corollary \ref{cor:rw_np2_nc}.
 \item $H$ is nonabelian, all its weights are real-valued, and the dimension of the nilpotent part is one. This case is treated in Proposition \ref{prop:rw_np1}, which contributes the cases (2.j)--(2.n). Note that case (1) from the Proposition corresponds to the cases (2.j) and (2.k) in Theorem \ref{thm:main}. This time, the admissible groups containing the connected ones are given in Proposition \ref{prop:rw_np1}. 
\end{itemize}

We now start with the classification. 
For subspaces $\h$ and $\k$ of $\mathfrak{gl}(d,\R)$, define 
$$
[\h,\k] = \sp \{ [X,Y] : X \in \h, Y \in \k \}.
$$
For any Lie subalgebra $\h$ of $\mathfrak{gl}(d,\R)$,  put $\h_{(0)} = \h^{(0} = \h$ and define $\h_{(k)}, \h^{(k)}$, $k = 1, 2, 3, \dots$ recursively by $\h_{(k)} = [\h_{(k-1)},\h_{(k-1)}]$ and $\h^{(k)} = [\h, \h^{(k-1)}$. Recall that $\h$ is said to be {\rm solvable} (resp. {\it nilpotent})  if there is some positive integer $n$ such that $\h_{(n)} = \{0\}$ (resp. $\h^{(n)} = \{0\}$. 

\

\



\begin{thm} {\rm (Lie's Theorem)} Any linear action of of a solvable Lie algebra on a complex vector space has a characteristic function.  






\end{thm}

The following consequence is almost immediate.

\begin{corollary}  \label{upper triangular} Let $\h$ be a solvable Lie subalgebra of $\gl(d,\R)$. 
Then $\h$ is conjugate with a Lie subalgebra of $\gl(d, \C)$ that consists of upper triangular matrices. Moreover, the set of weights for the standard representation of $\h$ is precisely the set of forms occurring on the diagonal of an upper triangular conjugate of $\h$.

\end{corollary}



Now suppose that $\h$ is a solvable Lie subalgebra of $\gl(3,\R)$, and let $W(\h)$ denote the set of weights for the standard representation of $\h$ . 
It is instructive to first treat the case where $\h$ is commutative. The following standard result shows that in this case every weight is in fact a characteristic function.


\begin{lemma} Let $\h$ be a commutative Lie subalgebra of $\mathfrak{gl}(3,\R)$.

\medskip

{\rm (a)} Given $\l \in W(\h)$, 
the subspace
$$
\Z_\l = \{ v \in \C^3 : \forall A \in \h, \exists p \in \N, (A^T-\l(A)I)^p v = 0\} 
$$
of $\C^3$ is $\h^T$-invariant, and there is a basis for $\Z_\l$ for which the restriction of $A^T$ to $\Z_\l$ has the form $\l I + N$ where $N$ is strictly lower-triangular.

\medskip

{\rm (b)} $\C^3 = \oplus \Z_\l$ where the sum is taken over the set $W(\h)$ of (distinct) weights. 


\end{lemma}















\



\

The preceding leads quickly to a classification of connected, abelian, open orbit, admissible subgroups of $GL(3,\R)$. Let $H$ be such a group, with Lie algebra $\h$. 
First suppose that there is a weight $\l$ that is not real-valued. Then $W(\h) =  \{ \gamma, \l, \overline{\l}\}$ where $\gamma$ is a real-valued form. Recall that $\Z_{\overline{\l}} = \overline{\Z}_\l$, and $\dim \Z_\l = \dim \Z_\gamma = 1$, so that $\h$ is simultaneously diagonalizable in $\C^3$. Write $\a = \Re \l$, $\b= \Im \l$. Observe that the restriction of the transpose action to $\Z_\l +  \overline{\Z}_\l$ has open orbits. Hence $\a$ and $\b$ are $\R$-linearly independent and $\a, \b , \gamma$ is a basis for the real linear dual of $\h$. Choose a basis for $\h$ dual to this basis, and choose a basis $e_1, e_2$, and $e_3$ for $\R^3$ such that $\Z_\gamma = \C e_1$ and  $\Z_\l = \C (e_2 +i e_3)$. Thus
$\h$ is conjugate to the span of 
$$
\left[\begin{matrix} 0&0&0\\0&1&0\\0&0&1\end{matrix}\right], \left[\begin{matrix} 0&0&0\\0&0&-1\\0&1&0\end{matrix}\right], \left[\begin{matrix} 1&0&0\\ 0&0&0\\ 0&0&0 \end{matrix}\right]
$$
and 
$H$ is conjugate to 
\[  H' = \left\{ \left[ \begin{array}{ccc} a & 0 & 0 \\ 0 & b \cos t & -b\sin t \\ 0 & b\sin t & b\cos t \end{array} \right] : a>0,b>0, t \in \mathbb{R}
 \right\}  .\]

Note that a
covering group of $H'$ for which there is one co-null open orbit is then $H' F$ where 
$$
F = \left\{ \left[\begin{matrix} \pm 1&0&0\\ 0&0&0\\ 0&0&0 \end{matrix}\right]\right\}.
$$
Now suppose that $W(\h)$ consists only of real-valued forms and let $n = |W(\h)|$. If $n = 3$ then $H$ is conjugate with  
\[ H' =  \left\{ \left[ \begin{array}{ccc} a_1 & 0 & 0 \\ 0 & a_2 & 0 \\ 0 & 0 & a_3 \end{array} \right] : a_i > 0 
\right\}  \]
and a covering group of $H'$ with one open orbit is $H'F$ where 
$$
F = \left\{ \left[\begin{matrix} \epsilon_1&0&0\\ 0&\epsilon_2&0\\ 0&0&\epsilon_3 \end{matrix}\right] : \epsilon_i = \pm 1
\right\}.
$$
If $n = 2$, then $\h$ conjugate to the span $\h'$ of 
$$
\left[\begin{matrix} 1&0&0\\0&1&0\\0&0&0\end{matrix}\right], \left[\begin{matrix} 0&1&0\\0&0&0\\0&0&0\end{matrix}\right], \left[\begin{matrix} 0&0&0\\ 0&0&0\\ 0&0&1 \end{matrix}\right],
$$
and the corresponding connected subgroup of $GL(3,\R)$ is 
\[ H' = \left\{ \left[ \begin{array}{ccc} a & b & 0 \\ 0 & a & 0 \\ 0 & 0 & c \end{array} \right] : a> 0,c> 0, b \in \mathbb{R},
\right\}  .\]
A covering group for which there is one co-null open orbit is then $H' F$ where 
$$
F = \left\{ \left[\begin{matrix}\epsilon_1&0&0\\ 0&\epsilon_1&0\\ 0&0&\epsilon_2 \end{matrix}\right] : \epsilon_i = \pm  1\right\}.
$$
Finally, suppose that there is exactly one real-valued weight $\l$. 
Choose a basis $e_1,e_2,e_3$ for $\R^3$ such that $\R e_3$ and $ \V_2 = \R e_2+ \R e_3$ are $\h^T$-invariant. Here $\n$ is codimension one in $\h$, and $H^T$ acts in $\R^3 / \R e_3$ with open orbits. Hence there must be $X \in \n$ such that $X^Te_1 = e_2 + c e_3$. Now  $Xe_2$ may or may not be zero, but for $Y \in \n$ such that $Y^T e_1 \in \Re_3$, $XY = YX$ implies that $Y^T e_2 = 0$. It follows that $ \dim \n = 2$, and if 
 $Xe_2 = 0$ also, then it is immediate that $\h$ is conjugate to the Lie algebra $\h'$ spanned by the elements $I, X, Y$ where 
$$
X^T = \left[\begin{matrix} 0&0&0\\1&0&0\\0&0&0\end{matrix}\right], Y^T = \left[\begin{matrix} 0&0&0\\ 0&0&0\\ 1&0&0 \end{matrix}\right],
$$
and $H$ is conjugate to 
\[ H' = \left\{ \left[ \begin{array}{ccc} a & b & c \\ 0 & a & 0 \\ 0 & 0 & a \end{array} \right] : a > 0, b, c \in \mathbb{R} \right\}  .\]
If $X^Te_2 = \delta e_3$ with $\delta  \ne 0$, then replacing $e_3 $ by $e_3' = \delta e_3$, we find that $\h$ is conjugate to $\h'$ spanned by $I, X,$ and $Y$ where 
$$
X^T=  \left[\begin{matrix} 0&0&0\\1&0&0\\0&1&0\end{matrix}\right], Y^T =  \left[\begin{matrix} 0&0&0\\ 0&0&0\\ 1&0&0 \end{matrix}\right],
$$
and $H$ is conjugate to the 
corresponding connected subgroup of $ GL(3,\R)$ given by
\[ H'= \left\{ \left[ \begin{array}{ccc} a & b & c \\ 0 & a & b \\ 0 & 0 & a \end{array} \right] : a > 0, b, c \in \mathbb{R} \right\}  .\]

\medskip
\noindent
In both cases, a covering group of $H'$ for which there is a single open orbit is $H' F$ where $F = \pm I$. Thus we have proved the following; see also \cite[Remark 20]{Fu98}.


\begin{proposition}\label{prop:ab} Let $H$ be a closed, connected abelian admissible subgroup of $GL(3,\R)$ for which $H^T$ has open orbits. Then $H$ is conjugate to exactly one of the following subgroups.

\begin{enumerate}
\item The diagonal group:
\[ H = \left\{ \left[ \begin{array}{ccc} a_1 & 0 & 0 \\ 0 & a_2 & 0 \\ 0 & 0 & a_3 \end{array} \right] : a_i >0 
\right\}  .\]

\item The block diagonal group $H = (\mathbb{R}^+ \cdot {\rm SO}(2)) \times \mathbb{R_+^*}$, that is
\[ H = \left\{ \left[ \begin{array}{ccc} a & 0 & 0 \\ 0 & b \cos t & -b\sin t \\ 0 & b\sin t & b\cos t \end{array} \right] : a> 0,b> 0, t \in \mathbb{R}
 \right\}  .\]

\item The group 
\[ H = \left\{ \left[ \begin{array}{ccc} a & b & 0 \\ 0 & a & 0 \\ 0 & 0 & c \end{array} \right] : a>  0,c>  0, b \in \mathbb{R}
\right\}  .\]

\item The group 
\[ H = \left\{ \left[ \begin{array}{ccc} a & b & c \\ 0 & a & b \\ 0 & 0 & a \end{array} \right] : a>  0 , b, c \in \mathbb{R}
\right\}  .\]

\item The group 
\[ H = \left\{ \left[ \begin{array}{ccc} a & b & c \\ 0 & a & 0 \\ 0 & 0 & a \end{array} \right] : a>0, b, c \in \mathbb{R} \right\}  .\]
\end{enumerate}

\end{proposition}

Of special interest are the irreducibly admissible groups. For abelian case, we have the following.

\begin{proposition} \label{prop:ab_nc} Let $H$ be an abelian irreducibly admissible subgroup of $GL(3.\R)$. Then $H$ is conjugate to exactly one of the following subgroups.

\begin{enumerate}
\item The diagonal group:
\[ H = \left\{ \left[ \begin{array}{ccc} a_1 & 0 & 0 \\ 0 & a_2 & 0 \\ 0 & 0 & a_3 \end{array} \right] : a_i \ne 0 
\right\}  .\]

\item The block diagonal group $H = (\mathbb{R}^+ \cdot {\rm SO}(2)) \times \mathbb{R_+^*}$, that is
\[ H = \left\{ \left[ \begin{array}{ccc} a & 0 & 0 \\ 0 & b \cos t & -b\sin t \\ 0 & b\sin t & b\cos t \end{array} \right] : a\ne 0,b\ne 0, t \in \mathbb{R}\right\}  .\]

\item The group 
\[ H = \left\{ \left[ \begin{array}{ccc} a & b & 0 \\ 0 & a & 0 \\ 0 & 0 & c \end{array} \right] : a\ne  0,c\ne  0, b \in \mathbb{R}
\right\}  .\]

\item The group 
\[ H = \left\{ \left[ \begin{array}{ccc} a & b & c \\ 0 & a & b \\ 0 & 0 & a \end{array} \right] : a\ne  0 , b, c \in \mathbb{R}
\right\}  .\]

\item The group 
\[ H = \left\{ \left[ \begin{array}{ccc} a & b & c \\ 0 & a & 0 \\ 0 & 0 & a \end{array} \right] : a \ne 0, b, c \in \mathbb{R},
c \not= 0 \right\}  .\]
\end{enumerate}
\end{proposition}


We now turn to the case where $H$ is non-abelian. Let $H$ be a closed, admissible, solvable, open orbit subgroup of $GL(3,\R)$ that is not abelian and let $\h$ be its Lie algebra. Then $\h$ is not commutative, and we have observed that $\dim \h = 3$ or $4$. 



\


Recall that if $\l \in W(\h)$ then $-\l$ is a weight for the transpose representation. We will be especially interested in those $\l \in W(\h)$ for which $-\l$ is a characteristic function for the transpose representation; denote the set of these weights by $C(\h)$. Thus for $\l \in C(\h)$, we have $f \in \C^3$ such that 
$$
A^T f = \l(A) f, \ \ \ A \in \h
$$
and we let $\E_\l$ denote the subspace of all $f$ satisfying the preceding. 
The {\it characteristic space} for $\h$ is defined by 
$$
\E(\h) = \oplus_{\l \in C(\h)} \E_\l;
$$
let $d(\h) $ denote the dimension (over $\C$) of $\E(\h)$. If $\l \in C(\h)$, then $\overline{\l} \in C(\h)$ also and $\E_{\overline{\l} } = \overline{\E}_\l$. It follows that $\E(\h)$ is the complexification of $\E(\h) \cap \R^3$. 

\

Denote by $c(\h)$ (resp. $w(\h) $) the dimension of the span of the characteristic functions (resp. weights) for the transpose action of $\h$.
If $\h'$ is conjugate to $\h$, then $c(\h') = c(\h), w(\h') = w(\h),$ and $d(\h') = d(\h)$. Since $\h$ is not commutative and open orbit, then $1 \le d(\h) \le 2$, and $1 \le c(\h)  \le w(\h) \le 2.$ 

\

Let $\n \subset \h $ be the subspace consisting of those elements for which all of the weights are zero. Since $[\h,\h] \subseteq \n$, then $\n$ is an ideal in $\h$, and since we are assuming that $\h$ is non-commutative, then $\n \ne \{0\}$. Note that $\dim \n = 3 - w(\h)$.

\

Let $\l\in C(\h)$ and suppose that $\l$ is not real-valued. Choose $f \in \E_\l$. The subspace spanned by $f, \overline{f}$ in $\C^3$ is the complexification of a two-dimensional subspace $\Z$ of $\R^3$ that is $\h^T$-invariant. Since $\R^3 / \Z$ is one-dimensional, then it is immediate that we have an ordered basis $(e, \overline{f}, f)$ for $\C^3$ with respect to which $A^T$ is given by the matrix 
\begin{equation}\label{matrix formula1}
\left[\begin{matrix} 
\gamma(A) &0&0\\ \overline{\tau}(A)  &\overline{\l}(A) &0\\    \tau(A)   &0 &\l(A)
\end{matrix}\right]
\end{equation}
for all $A \in \h$. Here $\gamma$ is a real linear form, and $\tau$ is a complex linear form on $\h$. Note that in this case $d(\h) \ge 2$.

%




\

\begin{proposition} \label{weight relations}  Let $\h$ be a non-commutative solvable, open orbit Lie algebra. 

\medskip
\noindent
{\rm (1)} 
If $C(\h)$ contains a real-valued form, then every form in $W(\h)$ is real-valued.

\medskip
\noindent
{\rm (2)} Suppose that $C(\h)$ contains a form $\l$ that is not real-valued, and write $\a = \Re \l, \b = \Im \l$. 
Then $\dim \n = 2$, and $\{ \a, \b, \gamma\}$ is $\R$-linearly dependent in $\h^*$. 

\end{proposition}

\begin{proof} To prove (1), suppose that $\gamma$ is a real-valued characteristic function and that there is a weight $\l$ that is not real-valued. As above, we have $e \in \R^3$, and $f, \overline{f} \in \C^3$, such that the matrix of $A^T$,  $A \in \h$ with respect to the ordered basis $(f, \overline{f}, e)$ is given by 
$$
  \left[\begin{matrix} \l(A)&0&0\\0&\overline{\l}(A) &0\\ \tau(A)&\overline{\tau}(A)&\gamma(A)\end{matrix}\right]
$$
where $\tau$ is a complex linear form on $\h$. Now since $\h^T$ is of open orbit type, then $\h^T$ acts with open orbits in the quotient $\R^3 / \R e$; hence $\l$ has real rank 2, that is, $\a$ and $\b$ are independent. Choose $B \in \h$ such that $\l(B) = i$, and observe that the transpose action of any element of $\n$ is determined by its action on $f$. Given a non-zero element $X\in \n$, $[B,X] \in \n$ also, but $ [B^T,X^T] f = (\gamma(B) - i)\tau(X) e$. Since $\gamma(B) - i \notin \R$, it follows that $\dim \n = 2$. This means that $\tau$ and $\overline{\tau}$ are linearly independent over $\R$. 

\

Let $X_1, X_2 \in \n$ such that $\tau(X_1) = 1, \tau(X_2) = i$. Write $f = e_1 + i e_2$ with $e_1, e_2 \in \R^3$. Now $(e_1, e_2, e)$ is an ordered basis for $\R^3$ and for $v \in \R^3$ write $v = v_1e_1 + v_2 e_2 + v_3 e$. If $v_1^2 + v_2^2 \ne 0$, then it is straightforward to check that $(v_1 X_2 - v_2 X_1) v = 0$, and hence that  $\exp \R (v_1 X_2 - v_2 X_1) $ is contained in the stabilizer of $v$. Thus almost all elements of $\R^3$ have a non-compact stabilizer, and $H$ is not admissible. Thus (1) is proved. 

\

To prove (2), choose an ordered basis $(e, \overline{f}, f )$ of $\C^3$ with $e \in \R^3$, such that  $A^T$ is given by 
$$
\left[\begin{matrix} 
\gamma(A) &0&0\\ \overline{\tau}(A)  &\overline{\l}(A) &0\\    \tau(A)   &0 &\l(A)
\end{matrix}\right]
$$
for all $A \in \h$. Choose $B \in \h$ such that $\l(B) \notin \R$. Since $\h$ is non-commutative, we have $X \in \n, X \ne 0$. Then $X^T e = \tau(X) f + \overline{\tau}(X) \overline{f}$, and since $X^T f = 0$, then $X^TB^T e = \gamma(B) X^T e$. Hence if $Y^T = [B^T,X^T] $, then $Y \in \n$ and 
$$
\begin{aligned}
Y^Te &= ( B^TX^T - X^T B^T)e 
= \l(B) \tau(X) f + \overline{\l}(B) \overline{\tau}(X) \overline{f} - \gamma(B) (\tau(X) f + \overline{\tau}(X) \overline{f}\\
&= (\l(B) - \gamma(B)) \tau(X) f + (\overline{\l}(B) - \gamma(B)) \overline{\tau}(X) \overline{f}
\end{aligned}
$$
so $\tau(Y) = (\l(B) - \gamma(B)) \tau(X)$. Since $\l(B) \notin \R$, then $\tau(Y)$ and $\tau(X)$ are linearly independent over $\R$, and hence $Y$ and $X$ are independent in $\n$. Thus $\dim \n \ge 2$, and since $\n$ is determined by the values of $\tau$ and $\overline{\tau}$, then $\dim \n = 2$.

Finally, since $\dim \h \le 4$ while $\dim \n \ge 2$, then $\dim (\text{span}_\R \{ \a, \b, \gamma\} )\le 2$. 
  \end{proof}

As a consequence of the preceding, we have the following corollary.

\begin{corollary} Let $\h$ be a non-commutative solvable, open orbit Lie subalgebra of $\gl(3,\R)$. 

\medskip

{\rm (1)} $\R^3$ contains a transpose-invariant subspace of codimension one. 

\medskip

{\rm (2)} Let $\V$ be a transpose-invariant subspace of codimension one, let $\l$ be the weight for the quotient action in $\R^3 / \V$, and let $\h_1 = \ker \l$. Then the transpose action of $\h_1$ in $\V$ is an open orbit action. 
\end{corollary} 
\

As in the case where $H$ is abelian, we describe the conjugacy classes of Lie algebras for solvable admissible $H$. First suppose that $C(\h)$ contains a form $\l$ that is not real-valued so that $W(\h) = \{\gamma, \l ,\overline{\l}\}$ where $\gamma$ is real-valued. By Lemma \ref{weight relations}, $\l$ and $\overline{\l}$ are characteristic functions, and we have an ordered basis $(e, \overline{f}, f)$  for $\C^3$ and a complex linear form $\tau$ on $\h$, such that the matrix of  $A^T$ with respect to the basis $(e, \overline{f}, f)$ is given by the form (\ref{matrix formula1}) 
for all $A \in \h$. 
Clearly we can take $e \in \R^3$, and writing $e = e_1, f = e_2 + i e_3$ with $e_2, e_3 \in \R^3$, then $(e_1 , e_2, e_3)$ is a basis for $\R^3$. Recall that also by Lemma \ref{weight relations}, $\dim(\text{span}_\R\{ \a, \b ,\gamma\}) = 1$ or $2$. 
 As before put $\a = \Re \l, \b = \Im \l$. Since $\l$ is not real-valued, then $\b \ne 0$. Hence there are three possibilities: (1) both $\a$ and $\gamma$ are multiples of $\b$, (2) $\a$ is a multiple of $\b$ but $\gamma$ is not, and (3) $\gamma$ is  multiple of $\b$ but $\a$ is not. Observe that these cases are ``conjugate invariant'': if $\h'$ is a solvable Lie algebra in $\gl(3,\R)$ that is conjugate to $\h$, then both $\h$ and $\h'$ fall within the same case. 
 

\


\noindent
(1) Suppose that both $\a$ and $\gamma$ are multiples of $\b$. Put $\mu = \Re \tau, \nu = \Im \tau$. 
Since the dimension of the span of $\a, \b, \gamma, \mu,$ and $\nu$ is at most 3.
We have 
$$
3 \le \dim \h = \dim \text{span}\{ \a,\b, \gamma, \mu,\nu\} \le 3
$$
so 
$\dim \h = 3$, and $\mu$ and $\nu$ are linearly independent. 
Choose $X_1, X_2$ in $\n$ such that with respect to the basis $(e_1,e_2,e_3)$, 
$$
X_1^T = \left[\begin{matrix} 0&0&0\\ 1&0&0\\0&0&0\end{matrix}\right], \ \ X_2^T = \left[\begin{matrix} 0&0&0\\ 0&0&0\\1&0&0\end{matrix}\right].
$$
Choose $A \in \h$ such that $\mu(A) = \nu(A) = 0$ while $\l(A) \ne 0$; using the form (\ref{matrix formula1}) for $A^T$, we compute that 
$$
[A, X_1 + i X_2] = (\gamma(A)-\a(A)  + i \b(A)) (X_1 + i X_2).
$$
Since $[X_1, X_2] = 0$ we get 
$$
\Delta_H(\exp t A) = e^{2t (\a(A) - \gamma(A))}
$$
while $|\det (\exp t A) | = e^{t(2\a(A) + \gamma(A))}$.  Since $H$ is admissible, then $\Delta_H \ne \det $ on $H$ so  $\gamma \ne 0$. It follows that $\{ \gamma, \mu,\nu\}$ is a basis for $\h^*$ and we fix the dual (ordered) basis $\{ A, X_1, X_2\}$ for $\h$. With $\l(A) = a + ib$, we now have
$$
[A, X_1 + i X_2] = (1-a  + i b) (X_1 + i X_2)
$$
as well as 
$[X_1,X_2] = 0$, characterizing the Lie algebra. To sum up, we have that $\h$ is conjugate to the Lie algebra spanned by $\{ A, X_1, X_2\}$ where
$$
A^T = \left[\begin{matrix} 1&0&0\\ 0&a &-b\\0&b &a \end{matrix}\right] , \ \ X_1^T = \left[\begin{matrix} 0&0&0\\ 1&0&0\\0&0&0\end{matrix}\right], \ \ X_2^T = \left[\begin{matrix} 0&0&0\\ 0&0&0\\1&0&0\end{matrix}\right].
$$
 

\


\

(2) Suppose that $\a$ is a multiple of $\b$, while $\gamma $ and $ \b$ are independent. Choose $B, C \in \h$ such that $\gamma(B) = 0,  | \l(B)| = 1,  \gamma(C) = 1, \l(C) = 0$. Note that $\b(B) \ne 0$, and that $\n = \ker \gamma \cap \ker \l$, being codimension-two in $\h$, is non-trivial.  

\


Now by Lemma \ref{weight relations}, $\dim \n = 2$, and using now the real basis $(e_1, e_2, e_3)$ for $\R^3$, a basis $X_1, X_2$ for $\n$ is given by 
$$
 X_1^T = \left[\begin{matrix} 0&0&0\\ 1&0&0\\0&0&0\end{matrix}\right], \ \ X_2^T = \left[\begin{matrix} 0&0&0\\0&0&0\\1&0&0\end{matrix}\right].
$$
Modifying $B$ and $C$ by linear combinations of $X_1$ and $X_2$ if necessary, we can assume that 
$Be_1 = Ce_1 = 0$. We compute that 
$$
[B, X_1 + i X_2] = -\overline{\l} (X_1 + i X_2), \ \ \ [C, X_1 + i X_2] = X_1 + i X_2
$$
characterizing the Lie algebra. It follows that 
$$
\Delta_H(\exp t B) = e^{2t\a}, \ \ \Delta_H(\exp t C) = e^{2t}
$$
while $|\det (\exp tB)|  = e^{2t\a}, \det \exp tC = e^t$. Thus $H$ satisfies the modularity condition. 
Whether $\h$ is the Lie algebra of an admissible group now depends only upon the stabilizers, which are generated by the annihilators in $\h$:
let $v = v_1 e_1 + v_ 2 e_2 + v_3 e_3$. Then the annihilator $\h(v)$ for generic $v \in \R^3$ is 
$$
\h(v) = \R \left(B + \Re \left(\frac{\l(v_2 + i v_3)}{2v_1}\right) X_1 + \Im \left(\frac{\l(v_2 + i v_3)}{2v_1}\right) X_2\right)
$$
By direct computation we verify that $\exp \h(v)$ is compact if and only if $\l$ is purely imaginary, that is, $\a = 0$. To sum up, if $\a$ is a multiple of $\b$ and $\gamma$ is independent of $\b$, and furthermore $\h$ is the Lie algebra of an admissible group, then $\h$ is conjugate with the Lie algebra spanned by $A, B, X_1, X_2$  where 
$$
B^T = \left[\begin{matrix} 0&0&0\\ 0&0 &-1\\0&1&0 \end{matrix}\right] , \ \  C^T = \left[\begin{matrix} 1&0&0\\ 0&0&0\\0&0&0\end{matrix}\right], \ \ X_1^T = \left[\begin{matrix} 0&0&0\\ 1&0&0\\0&0&0\end{matrix}\right], \ \ X_2^T = \left[\begin{matrix} 0&0&0\\0&0&0\\1&0&0\end{matrix}\right]
$$

\

(3)  Here we suppose that $\a$ is not a multiple of $\b$, while $\gamma$ belongs to the real span of $\a$ and $\b$. We choose $A, B \in \h$ such that $\a(A) = 1, \b(A) = 0, \a(B) = 0, \b(B) = 1$. 


\


Until the end of this case, we use the complex matrices with respect to the basis $(e, f, \overline{f})$ for $\C^3$. 
Again $\dim \n = 2$ and $\dim \h = 4$. 
Let $X_1, X_2 \in \n$ be given by $\tau(X_1) = 1$ and $\tau(X_2) = i$. Using (\ref{matrix formula1}) we compute $[A^T,X^T_j] = (1 - \gamma(A))X^T_j$ and $[B^T,X^T_1] = X^T_2 - \gamma(B) X^T_1, [B^T,X^T_2] = -X^T_1 -\gamma(B) X^T_2$. This can be expressed in the complexification of $\h$ more simply: with $X = X_1 + i X_2$, we get 
$$
[A^T,X^T] =  (1 - \gamma(A))X^T, \ \ [B^T,X^T] = (-\gamma(B) - i) X^T.
$$

To compute $\h(v)$, we first observe that  $\h(\pi(v)) = \ker\gamma = \k + \R (\gamma(B) A - \gamma(A) B)$. Now 
$$
\gamma(B) A^T - \gamma(A) B^T =  \left[\begin{matrix} 0&0&0\\0&\gamma(B) - i \gamma(A) &0\\ 0&0&\gamma(B) + i \gamma(A) \end{matrix}\right]. 
$$
Let $v \in \R^3$ such that $v_1 (v_2 + i v_3) \ne 0$. Then we have $a(v), b(v) \in \R$ such that $\h(v) = \R(\gamma(B) A - \gamma(A) B -a(v)X_1 - b(v)X_2)$ and we compute that 
$$
a(v) = \Re \left( \frac{(\gamma(B) - i \gamma(A)) (v_2 + i v_3)}{2v_1}\right) , \ \  b(v) = \Im \left( \frac{(\gamma(B) - i \gamma(A)) (v_2 + i v_3)}{2v_1}\right)
$$
This generates a torus if and only if $\gamma(B) = 0$. Coming back to the real basis $(e_1, e_2, e_3)$, we conclude that if the weights of the transpose action of $\h$ satisfy the conditions of 2.2, then there is a basis $(e_1, e_2, e_3)$ for $\R^3$ such that $\h$ is spanned by elements $A, B, X_1, X_2$, which are given with respect to this basis by
$$
A = \left[\begin{matrix} \gamma(A)&0&0\\ 0&1&0\\ 0&0&1\end{matrix}\right], B = \left[\begin{matrix} \gamma(B)&0&0\\ 0&0&-1\\0&1&0\end{matrix}\right], X_1 = \left[\begin{matrix} 0&0&0\\ 1&0&0\\0&0&0\end{matrix}\right], \ \ X_2 = \left[\begin{matrix} 0&0&0\\0&0&0\\1&0&0\end{matrix}\right]
$$
and $\h$ is admissible if and only if $\gamma(B) = 0$. Thus $\h$ is a span of elements $A, B, X_1, X_2$ and we have a basis for $\R^3$ such that these elements are given by 
$$
A^T = \left[\begin{matrix} a&0&0\\ 0&1&0\\ 0&0&1\end{matrix}\right], B^T = \left[\begin{matrix} 0&0&0\\ 0&0&-1\\0&1&0\end{matrix}\right], X_1^T = \left[\begin{matrix} 0&0&0\\ 1&0&0\\0&0&0\end{matrix}\right], \ \ X_2^T = \left[\begin{matrix} 0&0&0\\0&0&0\\1&0&0\end{matrix}\right]
$$
where $a\in \R, a \ne 0$. 

\

We summarize the preceding analysis. 

\

\begin{proposition} \label{prop:nrw} Let $\h$ be the Lie algebra of a connected, non-abelian admissible, open orbit subgroup $H$ of $GL(3,\R)$. Assume that $\h$ is solvable and has a weight that is not real valued. Then $W(\h) = \{ \l, \overline{\l}, \gamma\}$ where $\gamma $ is real valued, and both $\l$ and $\overline{\l}$ are characteristic functions. Moreover, exactly one of the following holds. 

\medskip
\noindent
{\rm (1)} $\l$ has real rank one, and $\gamma $ is a linear combination of $\l$ and $ \overline{\l}$. Then there is a unique $a+ i b \in \C \setminus \R$ such that $\h$ is conjugate with the Lie algebra spanned by 
$$
\left[\begin{matrix} 1&0&0\\ 0&a &-b\\0&b &a \end{matrix}\right] , \ \  \left[\begin{matrix} 0&1&0\\ 0&0&0\\0&0&0\end{matrix}\right], \ \ \left[\begin{matrix} 0&0&1\\ 0&0&0\\0&0&0\end{matrix}\right].
$$
The group $H$ is given by
$$
H = \left\{ \left[ \begin{matrix} e^s&t_1&t_2\\0&e^{as}\cos(bs)&-e^{as}\sin(bs)\\0&e^{as} \sin(bs) & e^{as} \cos(bs) \end{matrix} \right] : s, t_1, t_2 \in \R\right\}.
$$

\medskip
\noindent
{\rm (2)} $\l$ has real rank one, and $\gamma$ is not a linear combination of $\l$ and $ \overline{\l}$. Then $\h$ is conjugate to the Lie algebra spanned by 
$$
 \left[\begin{matrix} 0&0&0\\ 0&0 &-1\\0&1&0 \end{matrix}\right] , \ \   \left[\begin{matrix} 1&0&0\\ 0&0&0\\0&0&0\end{matrix}\right], \ \  \left[\begin{matrix} 0&1&0\\ 0&0&0\\0&0&0\end{matrix}\right], \ \  \left[\begin{matrix} 0&0&1\\0&0&0\\ 0&0&0\end{matrix}\right]
$$
Here we have
$$
H = \left\{ \left[ \begin{matrix} e^s&t_1&t_2\\0&\cos(r)&-\sin(r)\\0& \sin(r) & \cos(r) \end{matrix} \right] : r, s, t_1, t_2 \in \R\right\}.
$$

\medskip
\noindent
{\rm (3)} $\l$ has real rank two. Then there is a unique non-zero real number $a$ such that $\h$ is conjugate with the Lie algebra spanned by 
$$
A^T = \left[\begin{matrix} a&0&0\\ 0&1&0\\ 0&0&1\end{matrix}\right], B^T = \left[\begin{matrix} 0&0&0\\ 0&0&-1\\0&1&0\end{matrix}\right], X_1^T = \left[\begin{matrix} 0&1&0\\ 0&0&0\\0&0&0\end{matrix}\right], \ \ X_2^T = \left[\begin{matrix} 0&0&1\\0&0&0\\0&0&0\end{matrix}\right].
$$
In this case
$$
H = \left\{ \left[ \begin{matrix} e^{as}&t_1&t_2\\0&e^{s}\cos(r)&-e^{s}\sin(r)\\0&e^{s} \sin(r) & e^{s} \cos(r) \end{matrix} \right] : r, s, t_1, t_2 \in \R\right\}.
$$

\end{proposition}

In each of the above cases we see that $H$ has two open orbits, the orbits of $[\pm 1, 0, 0]^T$ in $\R^3$. Hence the following is immediate. 

\begin{corollary} \label{cor:nrw_nc}
 Let $\h$ be the Lie algebra of a connected, non-abelian admissible, open orbit subgroup $H$ of $GL(3,\R)$. Assume that $\h$ is solvable and has a weight that is not real valued. Let $E$ denote the identity matrix. Then $\{ -E, E \} H$ is irreduciablly admissible. 
 \end{corollary}


\

We now turn our attention to the second main case. Suppose now that $\h$ is the Lie algebra of a closed, solvable, admissible, open orbit subgroup of $GL(3,\R)$ and that every weight for the standard representation of $\h$ is real-valued. 
By Lemma \ref{upper triangular}, $\h$ is contained in the algebra of real upper-triangular $3\times 3$ matrices, so the matrix exponential is injective on $\h$. The image $H = \exp \h$ is the (unique) closed connected subgroup of $GL(3,R)$ corresponding to $\h$, and is an exponential Lie group. Now a one-parameter subgroup of an exponential group is non-compact, and for each $v \in \R^3$, the stabilizer $H(v)$ is connected. Since $H$ is open orbit admissible, then $\dim \h = 3$. 

\

Recall that $\E = \E(\h)$ is the span of all characteristic vectors for the transpose action; since $\h$ is not commutative, then $ \dim \E \le 2$. In the preceding analysis, we saw that if $\l \in W(\h)$ is not real-valued, then both $\l$ and $\overline{\l}$ belong to $C(\h)$, and hence $ \dim \E = 2$.  When all weights are real, then $\dim \mathcal E$ may be either 1 or 2, and this will lead a greater number of conjugacy classes.


\

As always, $\n = \cap_{\l \in W} \ker \l $ denotes the ideal of nilpotent elements in $\h$. Since $\h$ is an open orbit Lie algebra, then $\n \ne \h$; indeed $\dim \n = 1$ or $2$.  Another reason that the complex weight case is simpler is that by Lemma \ref{weight relations}, the presence of a weight that is not real implies that $\dim \n = 2$. In what follows we will distinguish conjugacy classes by these two characteristics: the dimension of $\E$, and the dimension of $\n$. It will also be helpful to identify the isomorphism classes of  three-dimensional solvable Lie algebras that are represented here. In addition to the Heisenberg Lie algebra, several non-nilpotent Lie algebras will occur. They are defined by writing the non-vanishing Lie brackets between basis elements. 

\

$\s_0 = \text{span}\{ A, C, Y\}$ with $ [A,Y] = Y$.

\

$\s_{1,c} = \text{span}\{ A, X, Y\}$ with $[A,X] = cX, [A,Y] = Y$. Here $0 <|c| \le 1$. 

\

$\s_2 = \text{span}\{ A, X, Y\}$ with $ [A,X] = X+Y, [A,Y] = Y$. 

\

In each case the connected group $H = \exp \h$ is easily obtained via the exponential map. An irreducibly admissible group is a finite extension of $H$. 

\

Suppose first that $\dim \n = 2$. Since $H$ is an open orbit group, there is a non-zero weight, of which every weight is a multiple. We describe the conjugacy classes of Lie algebras satisfying these requirements by identifying a basis of $\n$ and an element $A \notin \n$.


\

(1) Suppose that $\dim \E= 2$. Since $\h$ is of open orbit type then the quotient of the transpose action of $\h$ on the one-dimensional space $\R^3 / \E$ is given by a non-zero weight. 
Moreover, $\n$ acts trivially on $\E$. It follows immediately from these observations that $\h$ is conjugate to a Lie algebra spanned by $A, X, Y$ where $A$ is diagonal, and where $X$ and $Y$ are given by 
$$
X^T = \left[\begin{matrix}0&0&0\\ 1&0&0\\ 0&0&0\end{matrix}\right],  \ \ \ Y^T = \left[\begin{matrix}0&0&0\\ 0&0&0\\ 1&0&0\end{matrix}\right].
$$

\


Since $\h$ is non-commutative, then there are at least two distinct weights, that is, $A$ cannot be a multiple of the identity. If the weight for $\R^3 / \mathcal E$ coincides with a characteristic function,  
then 
$\h$ belongs to the class $\s_0$; indeed, $\h$ is conjugate to the Lie algebra spanned by  
 basis $A, X, Y$ where
$$
A^T =\left[\begin{matrix} \l&0&0\\ 0&\l &0\\ 0&0 &\l-1 \end{matrix}\right] 
$$
and where $\l \ne 0$. Note that here $C = X$ is central and $[A,Y] = Y$.

\

If $\dim \E = 2$ and there is a weight that is not a characteristic function, then 
$\h$ is conjugate to the Lie algebra spanned by  
$A, X, Y$ where
$$
A^T =\left[\begin{matrix} \l&0&0\\ 0&\l-1 &0\\ 0&0 &\l-b \end{matrix}\right] , \ \ X^T = \left[\begin{matrix}0&0&0\\ 1&0&0\\ 0&0&0\end{matrix}\right],  \ \ \ Y^T = \left[\begin{matrix}0&0&0\\ 0&0&0\\ 1&0&0\end{matrix}\right]
$$
where $\l \ne 0$ and $|b| \ge 1$. Thus $[A,X] = X, [A,Y]= bY$.


\

\

(2) Suppose that $\dim \n = 2$, and that $\dim \E = 1, \E = \R e_3$. Since $\E$ is invariant under the transpose representation, we have a representation $\a$ of $\h$ in $\R^3 / \E$; let $\k$ be the kernel of this quotient representation. Since the quotient map is open, then $\a$ is an open orbit action, and $\dim \k = 1$. Since $\h$ is open orbit, then $\n + \k \ne \h$, and hence $\k \subset \n$. Choose $X \in \n \setminus \k$ and $Y \in \n \cap \k, Y \ne 0$. The subspace $\V = X \R^3 + \E$ is invariant. With $e_1 \in \R^3$ such that $e_2 = X^Te_1 \notin \E$, then $\V = \text{span}\{e_2, e_3\}$ and it is easily seen that every element of $\h$ is upper triangular with respect to the ordered basis $(e_1, e_2, e_3)$. We can then modify this basis $\{X, Y\}$ of $\n$ so that they are given with respect to the basis $(e_1, e_2, e_3) $ by 
$$
X^T = \left[\begin{matrix}0&0&0\\ 1&0&0\\ 0&\delta&0\end{matrix}\right],  \ \ \ Y^T = \left[\begin{matrix}0&0&0\\ 0&0&0\\ 1&0&0\end{matrix}\right]
$$
where $\delta = 0 $ or $1$. 
We then have $A \in \h$ given by 
$$
A^T =\left[\begin{matrix} 1&0&0\\ 0&c_2&0\\ 0&\delta(A) &c_3 \end{matrix}\right].
$$

\

\


Now compute 
$$
 [A^T,X^T] = \left[\begin{matrix} 0&0&0\\ c_2 - 1&0&0\\ \delta(A)&(c_3 - c_2)\delta&0\end{matrix}\right]
$$
and since $[A^T,X^T] $ belongs to $\n^T$, it follows that 
\begin{equation}\label{brackets1} 
[A^T,X^T] = (c_2 - 1)X^T + \delta(A) Y^T, \ \ \ \ \  [A^T,Y^T] = (c_3-1)Y^T,
\end{equation}
and that 
\begin{equation}\label{relation1}
(c_3 - c_2)\delta = (c_2-1)\delta.
\end{equation}
The constants $c_2$ and $c_3$ express the relations among the weights of the transpose action of $\h$, and are therefore conjugate-invariant. Let $n$ be the number of distinct weights. In light of (\ref{brackets1}) and (\ref{relation1}), we observe the following.

\

\begin{itemize}

\item $n = 1$ if and only if  $c_2 = c_3 = 1$, whence $\delta(A) \ne 0$, and $\h$ is a Heisenberg Lie algebra. 

\item $n = 2$ and there are two weights for the quotient action in $\R^3 / \E$, if and only if $c_2 = c_3 \ne 1$. Here $\delta = 0$ and since $\V \ne \E$, then $\delta(A) \ne 0$. 

\item If $n=2$ and $c_2 \ne c_3$, then either $c_2 = 1$ or $c_3 = 1$, and hence $\delta =0$. But $c_2 \ne c_3$ also implies that $A$ is diagonalizable, whence $\dim \E = 2$, which is not the case here. 

\end{itemize}

\

The preceding observations lead to several isomorphism classes of Lie algebra. 

\

(a) Suppose that $n = 1$, that is, $c_2 = c_3 = 1$. As observed $\h$ is a Heisenberg Lie algebra and there are two one-parameter families of conjugacy classes. If $\delta = 0$, that is $X^2 = 0$, then  $\h$ is conjugate with the Lie algebra spanned by $A, X, Y$, given by
$$
A^T = \left[\begin{matrix} \l &0&0\\ 0&\l&0\\ 0&1&\l \end{matrix}\right], \ \ X^T =   \left[\begin{matrix} 0&0&0\\ 1&0&0\\ 0&0&0\end{matrix}\right], \ \ Y^T =   \left[\begin{matrix} 0&0&0\\ 0&0 &0\\ 1&0&0 \end{matrix}\right]
$$
where $\l \ne 0$. 
If $\delta  =1$ , then $\h$ is conjugate with the Lie algebra spanned by $A, X, Y$ given with respect to this basis by 
$$
A^T = \left[\begin{matrix} \l&0&0\\ 0&\l&0\\ 0&1&\l \end{matrix}\right], \ \ X^T =   \left[\begin{matrix} 0&0&0\\ 1&0&0\\ 0&1&0\end{matrix}\right], \ \ Y^T =   \left[\begin{matrix} 0&0&0\\ 0&0 &0\\ 1&0&0 \end{matrix}\right]
$$
Here again $\l \ne 0$ but now $X^2 \ne 0$, so these two families of conjugacy classes are distinct.

\

(b) Suppose that $n = 2$, but $c_2 = c_3$. The relation (\ref{relation1}) shows that $\delta = 0$. Here there is a two-parameter family of conjugacy classes each of which has a basis $B, U, V$ of $\h$ such that 
$$
[B,U] = U + V, \ \ \ \ \ [B,V] = V.
$$
Namely $\h$ is conjugate with the Lie algebra spanned by $B, U, V$ given by
$$
B^T = \l A^T =  \left[\begin{matrix} \l&0&0\\ 0&\l+1&0\\ 0&\nu&\l+1 \end{matrix}\right], \ \ U^T =   \left[\begin{matrix} 0&0&0\\ 1&0&0\\ 0&0&0\end{matrix}\right], \ \ V^T =   \left[\begin{matrix} 0&0&0\\ 0&0 &0\\ \nu&0&0 \end{matrix}\right]
$$
where $\l = 1 / (c_2 - 1)$, and $\nu \ne 0$.

\

The third observation above shows that there is only one more case:




\medskip
\noindent
(c) $n = 3$, so that $c_2 \ne c_3$ and $c_2 \ne 1, c_3 \ne 1$. The basis $(e_2, e_3)$ for $ \V$ can be chosen so that $A$ is diagonal, and since $\dim \E = 1$, we must have $\delta = 1$. Putting $\l = 1 / (1- c_2 )$, then the relations (\ref{relation1}) show that $\h$ is conjugate with the Lie algebra spanned by $B, X, Y$ where
$$
B^T = \l A =  \left[\begin{matrix} \l&0&0\\ 0&\l-1&0\\ 0&0&\l-2 \end{matrix}\right], \ \ X^T =   \left[\begin{matrix} 0&0&0\\ 1&0&0\\ 0&1&0\end{matrix}\right], \ \ Y^T =   \left[\begin{matrix} 0&0&0\\ 0&0 &0\\ 1&0&0 \end{matrix}\right].
$$
and $[B,X] = X, [B,Y] = 2Y$.

\

We sum up as follows. 

\

\begin{proposition} \label{prop:rw_np2} Let $H$ be a closed, connected, admissible, open orbit subgroup of $GL(3,\R)$ whose Lie algebra $\h$ is non-commutative solvable with real weights. Suppose further that the dimension of the subalgebra of nilpotent matrices in $\h$ is two. Then $\h$ is conjugate to exactly one of the following Lie algebras. 

\medskip

{\rm (1)} If $\dim \E = 2$ and every weight is also a characteristic function, then 
$\h$ is conjugate to the Lie algebra spanned by  
 basis $A, C, Y$ where
$$
A =\left[\begin{matrix} \l&0&0\\ 0&\l &0\\ 0&0 &\l-1 \end{matrix}\right], \ \ \ C = \left[\begin{matrix}0&1&0\\ 0&0&0\\ 0&0&0\end{matrix}\right],  \ \ \ Y^T = \left[\begin{matrix}0&0&1\\ 0&0&0\\ 0&0&0\end{matrix}\right].
$$
and where $\l \ne 0$. Here $C$ is central and $[A,Y] = Y$, so $\h$ belongs to the isomorphism class $\s_0$. Computing the matrix exponential one finds 
$$
H = \exp \h = \left\{ \left[\begin{matrix} a^\l & t_1&t_2\\ 0& a^{\l}&0\\ 0&0&a^{\l-1}\end{matrix}\right] : a > 0, t_1, t_2 \in \R \right\}.
$$

\medskip

{\rm (2)}  If $\dim \E = 2$ and there is a weight that is not a characteristic function, then 
 $\h$ is conjugate to the Lie algebra spanned by  
$A, X, Y$ where
$$
A =\left[\begin{matrix} \l&0&0\\ 0&\l-c &0\\ 0&0 &\l-1 \end{matrix}\right] , \ \ X = \left[\begin{matrix}0&1&0\\ 0&0&0\\ 0&0&0\end{matrix}\right],  \ \ \ Y= \left[\begin{matrix}0&0&1\\ 0&0&0\\ 0&0&0\end{matrix}\right]
$$
where $\l \ne 0$ and $0<|c|\le1$. (Note that (1) is a limiting case of the above.) Thus $[A,X] = cX, [A,Y]= Y$, so $\h$ belongs to the isomorphism class $\s_{1,b}$.  Here one finds that $H = \exp \h$ is given by
$$
H = \left\{ \left[\begin{matrix} a^\l & t_1 &t_2\\ 0& a^{\l-c} &0\\ 0&0&a^{\l-1}\end{matrix}\right] : a > 0, t_1, t_2 \in \R\right\}.
$$


 \medskip
 
{\rm (3)}  If $\dim \E = 1$ and there is only one weight, 
 then  $\h$ is conjugate with the Lie algebra spanned by $A, X, Y$, given either by
$$
A = \left[\begin{matrix} \l &0&0\\ 0&\l&1\\ 0&0&\l \end{matrix}\right], \ \ X =   \left[\begin{matrix} 0&1&0\\ 0&0&0\\ 0&0&0\end{matrix}\right], \ \ Y =   \left[\begin{matrix} 0&0&1\\ 0&0 &0\\ 0&0&0 \end{matrix}\right]
$$
or 
$$
A = \left[\begin{matrix} \l&0&0\\ 0&\l&1\\ 0&0&\l \end{matrix}\right], \ \ X =   \left[\begin{matrix} 0&1&0\\ 0&0&1\\ 0&0&0\end{matrix}\right], \ \ Y=   \left[\begin{matrix} 0&0&1\\ 0&0 &0\\ 0&0&0 \end{matrix}\right]
$$
where $\l \ne 0$. In both cases $\h$ is a Heisenberg Lie algebra. In the first case
$$
H = \exp \h = \left\{ \left[\begin{matrix} a^\l& t_1&t_2  \\ 0& a^\l &a^\l \ln a \\ 0&0&a^\l \end{matrix}\right] : a >0, t_1, t_2 \in \R\right\}.
$$
while in the second case, 
$$
H = \exp \h = \left\{ \left[\begin{matrix} a^\l& t_1&t_2  \\ 0& a^\l &t_1 + a^\l \ln a \\ 0&0&a^\l \end{matrix}\right] : a >0, t_1, t_2 \in \R\right\}.
$$


 \medskip
 
{\rm (4)}   If $\dim \E = 1$ and there are two weights for the quotient action in $\R^3 / \E$, then $\h$ is conjugate with the Lie algebra spanned by 
$$
A =\left[\begin{matrix} \l&0&0\\ 0&\l-1&\delta \\ 0&0 &\l-1 \end{matrix}\right] , \ \ X = \left[\begin{matrix} 0&1&0\\ 0&0&0\\ 0&0&0\end{matrix}\right], \ \ Y = \left[\begin{matrix} 0&0&1\\ 0&0&0\\ 0 &0&0\end{matrix}\right].
$$
Here $\h$ belongs to the isomorphism class $\s_2$. In this case $H = \exp \h$ is given by 
$$
H = \left\{ \left[\begin{matrix} a^\l& t_1 &t_2\\ 0& a^{\l-1}&\delta  a^{\l-1}\ln a\\ 0&0&a^{\l-1}\end{matrix}\right] : a>0,  t_1, t_2 \in \R\right\}.
$$ where $\lambda, \delta \in \mathbb{R}^*$.

\medskip

{\rm (5)}   If $\dim \E = 1$ and there are three distinct weights, then $\h$ is conjugate with the Lie algebra spanned by $A, X, Y$ where
$$
A =  \left[\begin{matrix} \l&0&0\\ 0&\l-1&0\\ 0&0&\l-2 \end{matrix}\right], \ \ X =   \left[\begin{matrix} 0&1&0\\ 0&0&1\\ 0&0&0\end{matrix}\right], \ \ Y =   \left[\begin{matrix} 0&0&1\\ 0&0 &0\\ 0&0&0 \end{matrix}\right],
$$
with $\lambda \not= 0$. 
Here $\h$ belongs to the isomorphism class $\s_{1,1/2}$. Exponentiating and simplifying gives
$$
H = \exp \h = \left\{ \left[\begin{matrix} a^\l& t_1 &t_2\\ 0& a^{\l-1}& t_1 a^{\l-2}\\ 0&0&a^{\l-2}\end{matrix}\right] : a > 0, t_1, t_2 \in \R\right\}.
$$

\end{proposition}

\

The same argument as for Corollary \ref{cor:nrw_nc} is applicable here as well:
\begin{corollary} \label{cor:rw_np2_nc} Let $H$ be a closed, connected, admissible, open orbit subgroup of $GL(3,\R)$ whose Lie algebra $\h$ is non-commutative solvable with real weights. Suppose further that the dimension of the subalgebra of nilpotent matrices in $\h$ is two. Let $E$ denote the identity matrix in $GL(3,\R)$. Then $\{E, -E\} H$ is irreducibly admissible. 
\end{corollary}

Remark that the group $FH$ is also irreducibly admissible, where
$$
F = \left\{  \left[\begin{matrix} \epsilon &0&0\\ 0&1&0\\ 0&0&1 \end{matrix}\right] : \epsilon\in \{ -1, 1\} \right\}.
$$

\

Now we analyze the situation where $\dim \n = 1$. The following shows that if $\dim \n = 1$, then only one isomorphism class of Lie algebras occurs. 

\

\begin{lemma} \label{oneLiealgebra} Suppose that $\dim \n = 1$. Then $\h$ belongs to the class $\s_0$, and $[\h,\h] = \n$. 

\end{lemma} 

\begin{proof} Since $\h$ is non-commutative and 
$\dim \n = 1$, then $[\h,\h] = \n$. Write $\n = \R Y$; it is enough to show that $Y$ is not central in $\h$: since $\n = [\h,\h]$, then we have $A \in \h$ such that $[A,Y] = Y$, and if $B$ is independent of $A$ and $Y$, then $[A,B] = cY$ for some $c \in \R$ and $C = B - cY$ is central. 

\

To see that $Y$ is not central in $\h$, suppose the contrary; then $\h$ is Heisenberg. Let $\overline{\h}$ be the image of $\h$ in $\gl(2,\R)$ given by the quotient action of $\h^T$ on $\R^3 / \R e_3$. Since $\gl(2,\R)$ does not contain a Heisenberg algebra, then the image of $Y$ in $\overline{\h}$ is zero and $\overline{\h}$ is commutative. Now there are exactly two conjugacy classes of commutative Lie algebras in $\gl(2,\R)$ that are of open orbit type: 
$\overline{\h}$ is conjugate to the Lie algebra spanned either by 
$$
\left[ \begin{matrix} 1&0\\0&0\end{matrix}\right], \ \ \ \left[ \begin{matrix} 0&0\\0&1\end{matrix}\right]
$$
or by 
$$
\left[ \begin{matrix} 1&0\\0&1\end{matrix}\right], \ \ \ \left[ \begin{matrix} 0&0\\1&0\end{matrix}\right].
$$
In the first case $\h$ is conjugate to a Lie algebra spanned by elements of the form
$$
P = \left[ \begin{matrix} 1&0&0\\0&0&0\\ \nu_1&\delta_1&c\end{matrix}\right], \ \ \ Q =  \left[ \begin{matrix} 0&0&0\\0&1&0\\ \nu_2&\delta_2&d \end{matrix}\right] , \ \ \ Y =  \left[ \begin{matrix} 0&0&0\\0&0&0\\ \nu_3&\delta_3&0 \end{matrix}\right],
$$
where $Y$ commutes with $P$ and $Q$ and $[P,Q] \in \R Y$; a simple computation shows that this is impossible. A similar argument shows that the second case also leads to a contradiction.

\end{proof}

\

We continue to assume that $\n = [\h,\h] = \R Y$ for some $Y\ne 0$. Observe that $\E$ is $\h^T$-invariant and $Y^T \E = 0$. Recall that 
 $\dim \E \le 2$ since $\h$ is non-commutative. 

\

First suppose that $\dim \E =  2$, and that there are two linearly independent characteristic functions. Here we have an ordered basis $(e_1, e_2, e_3)$ for $\R^3$ and a basis $A, B, Y$ for $\h$ such that with respect to the basis $(e_1, e_2, e_3)$, 
$$
A^T = \left[\begin{matrix} a&0&0\\ \mu(A)&1 &0\\\nu(A)&0 &0 \end{matrix}\right] , \ \ B^T = \left[\begin{matrix} b&0&0\\ \mu(B)&0&0\\\nu(B)&0&1\end{matrix}\right], \ \ Y^T = \left[\begin{matrix} 0&0&0\\ \mu(Y)&0&0\\\nu(Y) &0&0\end{matrix}\right].
$$
The fact that  $[A,Y] \in \R Y$ and $[B,Y] \in \R Y$ implies that $\mu(Y)\nu(Y) = 0$. We may then assume that the basis $(e_1, e_2, e_3)$ and $A, B, Y$ for $\h$ are chosen so that 
$$
A^T = \left[\begin{matrix} a&0&0\\ 0&1 &0\\\nu(A)&0 &0 \end{matrix}\right] , \ \ B^T = \left[\begin{matrix} b&0&0\\ 0&0&0\\\nu(B)&0&1\end{matrix}\right], \ \ Y^T = \left[\begin{matrix} 0&0&0\\ 1&0&0\\ 0 &0&0\end{matrix}\right].
$$
Here $[A,Y] = (a-1) Y, [B,Y] = bY$; since $Y$ is not central, then either $a \ne 1 $ or $b \ne 0$ must hold. Since $[A,B] e_1 \in \R e_3$ but also $[A,B] \in \R Y$, then $[A, B] = 0$. Hence the center of $\h$ is $\R C$ where $C$ is given by 
$$
C^T = bA^T + (1-a) B^T =  \left[\begin{matrix} b&0&0\\ 0&b &0\\ \nu(C)&0 &1-a \end{matrix}\right].
$$
If $a+ b \ne 1$, then (again modifying the basis of $\R^3$) we may assume that $\nu(C) = 0$, and since $[A,C] = [B,C] = 0$, then $\nu(A) = \nu(B) = 0$ also.  If $a \ne 1$ then replacing $A$ by $\l A$ where  $\l = 1/(a-1)$ gives 
$$
A' =\left[\begin{matrix} \l+1&0&0\\ 0&\l&0\\ 0&0 &0 \end{matrix}\right] 
$$
while if $a = 1$ and $b \ne 0$, then we can take 
$$
A' = (1/b)B = \left[\begin{matrix} 1&0&0\\ 0&0&0\\ 0&0 &1/b\end{matrix}\right] .
$$
These cases are conjugate-invariant since the rank of $C$ is either 1 or 3 when $a \ne 1$, while if $a = 1$ and $b \ne 0$ then the rank of $C$ is two. 
Suppose now that $a + b = 1$; both $a \ne 1$ and $b \ne 0$ hold. Here the center is $\R C$ where 
$$
C^T = \left[\begin{matrix} 1&0&0\\ 0&1 &0\\ \nu(C) &0 &1 \end{matrix}\right] 
$$
for some $\nu(C) \in \R$. 
 If $a = 0$ and $b = 1$, then we get that $\h$ is conjugate to the Lie algebra spanned by $A, C, Y$ where 
 $$
 A^T = \left[\begin{matrix} 0&0&0\\ 0&1 &0\\  \nu(A) &0 &0 \end{matrix}\right] 
 $$
and $\nu(A) \in \R$. Observe that in this case, for $v \in \R^3$ such that $v_1 \ne 0$, the matrix $M(v)$ has rank 3 if and only if $\nu(A) \ne 0$. Thus the assumption that $\h$ is an open orbit subalgebra means that $\nu(A) \ne 0$.  If $a \ne 0$, then again replacing $A$ by $\l A$ where $l = 1 / (a-1)$ we get that  $\h$ is conjugate to the Lie algebra spanned by $A, C, Y$ where $C = I$ and 
$$
 A^T = \left[\begin{matrix} \l +1&0&0\\ 0&\l &0\\ 0 &0 &0 \end{matrix}\right] .
 $$
 Note that in this case both $\l$ and $\l+1$ are non-zero. 

\

Continuing the analysis when $\dim \n = 1$, suppose now that $\dim \E  = 2$, but that the characteristic functions are linearly dependent. Here again $Y \E =0$, and of course there are two linearly independent weights. Hence we can choose an ordered basis $(e_1, e_2, e_3)$ for $\R^3$ such that $e_2, e_3$ spans $\E$ and  $Ye_1 = e_2$. We then have a basis $A, B, Y$ for $\h$ given by 
$$
A^T =\left[\begin{matrix} 0&0&0\\ 0&1&0\\ \nu(A)&0 &a \end{matrix}\right] , \ \ B^T = \left[\begin{matrix} 1&0&0\\ 0&0&0\\ \nu(B)&0&0\end{matrix}\right], \ \ Y^T = \left[\begin{matrix} 0&0&0\\ 1&0&0\\ 0 &0&0\end{matrix}\right]
$$
with respect to the chosen basis. Replacing $e_1$ by $e_1 - \nu(A)e_3$ we may assume that $\nu(A) = 0$. One finds that $[A,Y] = -Y$ and $[B,Y] = Y$, and since $[A,B] \in \R Y$, then $[A,B] = 0$. So with $C = A+B$, we have $\R C$ is the center of $\h$. Note that the relation $[A,B] = 0$ implies that $\nu(A) + a \nu(B) = 0$. If $a \ne 1$, then $C$ is diagonalizable and since $A$ commutes with $C$, then $A$ is simultaneously diagonalizable; thus
 $\h$ is conjugate to the Lie algebra spanned by $A, C, Y$ where
$$
A^T =\left[\begin{matrix} 0&0&0\\ 0&1&0\\ 0&0 &a \end{matrix}\right] , \ \ C^T = \left[\begin{matrix} 1&0&0\\ 0&1&0\\ 0 &0&a\end{matrix}\right], \ \ Y^T = \left[\begin{matrix} 0&0&0\\ 1&0&0\\ 0 &0&0\end{matrix}\right].
$$
 If $a = 1$, then $\nu(A) + \nu(B) = 0$ and $C$ is the identity matrix. Thus
$\h$ is conjugate to the Lie algebra spanned by $A, C, Y$ where
$$
A^T =\left[\begin{matrix} 0&0&0\\ 0&1&0\\ 0&0 &1 \end{matrix}\right] , \ \ C^T = \left[\begin{matrix} 1&0&0\\ 0&1&0\\ 0&0&1\end{matrix}\right], \ \ Y^T = \left[\begin{matrix} 0&0&0\\ 1&0&0\\ 0 &0&0\end{matrix}\right].
$$

Next we consider the situation where $\dim \n = \dim E = 1 $; write $\E = \R e_3$ for some $e_3 \in \R^3$, and we consider the quotient action in $\R^3 / \E$. As before we consider the quotient of the transpose representation of $\h$ acting in $\R^3 / \E$, and let $\k$ denote its kernel. As before we have $\dim \k = 1$. 

\

\begin{lemma} \label{quotient is commutative} Let $\h$ be a solvable open orbit Lie subalgebra of $\gl(3,\R)$ and suppose that $\dim \n = \dim \E = 1$. Then the action of $\h$ in $\R^3 / \E$ is diagonalizable.  

\end{lemma}

\begin{proof} Since $\dim \n = 1 $ also, then either $\n = \k$, or $\n \cap \k = \{0\}$. If $\n = \k$, then $\h / \k$ is a commutative Lie aubalgebra of $\gl(2,\R)$, and this case there are only two possibilities for the quotient of the transpose action in $\R^3 / \E$:  either it has two independent characteristic functions (i.e., is diagonalizable), or just one weight.  

\

Now suppose that the lemma is false. Then either $\n = \k$ and the action of $\h$ in  $\R^3 / \E$  has only one weight, or $\n \ne \k$. We dispense with these cases one at a time. 

\

Suppose that $\n = \k$ and that the action of $\h$ in  $\R^3 / \E$  has only one weight. We have an ordered basis $(e_1, e_2, e_3)$ for $\R^3$ with respect to which each $A \in \h$ is given by  
$$
A^T =   \left[\begin{matrix} \l(A)&0&0\\ \mu(A)&\l(A) &0\\ \nu(A)&\delta(A)&\l_2(A)\end{matrix}\right]
$$
where again $\mu, \nu, $ and $\delta$ are real linear forms on $\h$. The open orbit assumption implies that $\l \ne 0$ and that $\mu$ is independent of $\l$. Hence we have $A, B, Y \in \h$ of the form
$$
A^T =\left[\begin{matrix} 1&0&0\\ 0&1&0\\ \nu(A)&\delta(A) &a \end{matrix}\right] , \ \ B^T = \left[\begin{matrix}0&0&0\\ 1&0&0\\ \nu(B)&\delta(B)&b\end{matrix}\right],  \ \ \ Y^T = \left[\begin{matrix}0&0&0\\ 0&0&0\\ \nu(Y)&\delta(Y)&0\end{matrix}\right],  
$$
Since $\dim \n = 1$ then $b \ne 0$. Now we compute that 
$$
[B^T,Y^T] =  \left[\begin{matrix}0&0&0\\ 0&0&0\\ b \nu(Y) - \delta(Y) & b \delta(Y)&0\end{matrix}\right],  
$$
and since $[B,Y] \in \R Y$ this leads to $\delta(Y) = 0$. We choose $Y$ so that $\nu(Y) = 1$, and $A$ and $B$ so that $\nu(A) = \nu(B) = 0$. Since $b \ne 0$, then the basis $(e_1, e_2, e_3)$ can be chosen so that $\delta(B) = 0$. Finally, observe that since $[A,B] \in \R Y$ then $\delta([A,B]) = 0$; but a direct computation shows that $\delta([A,B]) = -d \delta A$, so $\delta(A) = 0$. but now $e_2$ and $e_3$ are characteristic vectors for $\h$, a contradiction to our assumption that $\dim \E = 1$. 


\

Next assume that $\n \ne \k$, that is, that the quotient action of $\h$ on $\R^3 / \E$ is not commutative. Writing $\E = \R e_3$, Choose $C \in \k$ so that $C e_3 = e_3$. Since $\overline{Y} \ne 0$, we may choose an ordered basis  $(e_1, e_2, e_3) $ for $\R^3$ so that $\h$ is upper triangular, and so that $Y e_1 = e_2$. Now $[C,Y] \in \n \cap \k = \{0\}$ so $[C,Y] = 0$ and $YC = 0$, so $0 = Y Ce_1 = CYe_1 = Ce_2$. Similarly $[A,C]= 0$. Thus the center of $\h$ is $\R C$. It follows that both $e_2$ and $e_3$ are characteristic vectors for $\h$, again contradicting the assumption that $\dim \E = 1$.



\end{proof}

\

We proceed to describe the distinct conjugacy classes of $\h$ when $\dim \E = 1$ and $\dim \n = 1$. 
By Lemma \ref{quotient is commutative} we have an ordered basis $(e_1, e_2, e_3)$ for $\R^3$ with respect to which each $A \in \h$ is given by  
$$
A^T =   \left[\begin{matrix} \l_1(A)&0&0\\ 0&\l_2(A) &0\\ \nu(A)&\delta(A)&\l_3(A)\end{matrix}\right]
$$
where $\l_i, \mu, \nu, \delta$ are real linear forms on $\h$, and where the characteristic function $\l_3$ is dependent upon $\l_1$ and $\l_2$. 
Now choose $A, B \in \h$ such that $\l_1(A) = \l_2(B) = 1, \l_2(A) = \l_1(B) = 0$. Given $Y$ so that $\n = \R Y$, 
then a computation of $[A^T, Y^T] \in \R Y^T$ shows that either $\nu(Y) = 0$ or $\delta(Y)= 0$. Since the basis elements  $e_1, e_2$ are interchangeable, then we may assume that $\nu(Y) = 1, \delta(Y) = 0$, and then we may assume that $\nu(A) = \nu(B) = 0$. 
Thus we have a basis $(e_1, e_2, e_3)$ for $\R^3$ and $A, B, Y \in \h$ given with respect to this basis by
$$
A^T =\left[\begin{matrix} 1&0&0\\ 0&0&0\\ 0&\delta(A) &a \end{matrix}\right] , \ \ B^T = \left[\begin{matrix}0&0&0\\ 0&1&0\\ 0&\delta(B)&b\end{matrix}\right],  \ \ \ Y^T = \left[\begin{matrix}0&0&0\\ 0&0&0\\ 1&0&0\end{matrix}\right].
$$
We compute that $[A,Y] = (1-a) Y, [B,Y] =-bY$; 
since $Y$ is not central then either $a-1$ or $b$ is non-zero, 
and since $[A,B] \in \R Y$ while $[A,B] e_1 = 0$, then $[A,B] = 0$. Hence the center of $\h$ is spanned by the element $C =bA + (1-a) B$ and with respect to our basis for $\R^3$,  $C^T$ is given by 
\begin{equation}\label{central1}
C^T = \left[\begin{matrix} b&0&0\\ 0&1-a&0\\ 0&\delta(C) &b \end{matrix}\right].
\end{equation}
If $a+b \ne 1$, then $C$ is diagonalizable on the span of $e_2 $ and $e_3$; but then, $C$ being central in $\h$, both $e_2$ and $e_3$ would be characteristic vectors for all of $\h$. Since we assume now that $\dim \E = 1$, this is not the case. Hence $a+b = 1$, and we may assume that 
$$
C^T = \left[\begin{matrix} 1 &0&0\\ 0&1&0\\ 0&\delta(C)  &1 \end{matrix}\right].
$$
Now with $A^T$ as above, $[A,C] = 0$ implies $a \delta(C) = 0$. If $a \ne 0$, then $\delta(C) = 0$ and $C = I$, but also (modifying the element $e_2$) we get $\delta(A) = 0$, again contradicting our assumption that $\dim \E = 1$. Thus $a = 0$, and 
with respect to the basis $(e_1, e_2, e_3)$,  
$$
A^T =\left[\begin{matrix} 1&0&0\\ 0&0&0\\ 0&\delta(A) &0\end{matrix}\right], \ \ \ C^T = \left[\begin{matrix} 1 &0&0\\ 0&1&0\\ 0&\delta(C)  &1 \end{matrix}\right], \ \ \ Y^T = \left[\begin{matrix}0&0&0\\ 0&0&0\\ 1 &0&0\end{matrix}\right].
$$
Since $\dim \E = 1$, then either $\delta(A) \ne 0 $ or $\delta(C) \ne 0$.


\

\

We sum up the preceding. 

\

\begin{proposition} \label{prop:rw_np1} Let $H$ be a closed, connected, admissible, open orbit subgroup of $GL(3,\R)$ whose Lie algebra $\h$ is non-commutative solvable with real weights. Suppose further that the dimension of the subalgebra of nilpotent matrices in $\h$ is one. Then $\h$ is conjugate to a Lie algebra spanned by matrices $A, C$, and $Y$ where $[A,Y] = Y$, $\R C$ is the center of $\h$, and where the matrices $A, C$, and $Y$ are exactly one of the following. 

\medskip

{\rm (1)} If there are two linearly independent characteristic functions and three distinct weights, then there are two families of conjugacy classes. If the rank of a non-zero central element of $\h$ is odd, then 
 $$
A^T =\left[\begin{matrix} \l+1&0&0\\ 0&\l&0\\ 0&0 &0 \end{matrix}\right] , \ \ C^T = \left[\begin{matrix} c&0&0\\ 0&c&0\\ 0&0&1\end{matrix}\right], \ \ Y^T = \left[\begin{matrix} 0&0&0\\ 1&0&0\\ 0 &0&0\end{matrix}\right]
$$
while if the rank of a non-zero central element of $\h$ is 2, then 
$$
A^T =\left[\begin{matrix} 1&0&0\\ 0&0&0 \\ 0&0 &\l \end{matrix}\right] , \ \ C^T = \left[\begin{matrix} 1&0&0\\ 0&1&0\\ 0&0&0\end{matrix}\right], \ \ Y^T = \left[\begin{matrix} 0&0&0\\ 1&0&0\\ 0 &0&0\end{matrix}\right]
$$
where $\l \ne 0$, $c \ne 1$. In the first case the group $H=\exp \h$ is
$$
H = \left\{ \left[\begin{matrix} a^{\l+1} b^c & t&0\\ 0& a^\l b^c&0\\ 0&0&b \end{matrix}\right] : a>0, b>0, t \in \R\right\}
$$
while in the second, 
$$
H = \left\{ \left[\begin{matrix} a b & t&0\\ 0& b&0\\ 0&0&a^\l  \end{matrix}\right] : a>0, b>0, t \in \R\right\}
$$
An irreducibly admissible group is here given by  $F H$ where
$$
F = \left\{ \left[\begin{matrix} \epsilon_1 & 0&0\\ 0& 1 &0\\ 0&0&\epsilon_2\end{matrix}\right] : \epsilon_1, \epsilon_2 \in \{ -1, 1\}\right\}.
$$


\medskip

{\rm (2)} If there are two linearly independent characteristic functions and two distinct weights, then again there are two families of conjugacy classes. 
 $$
A^T =\left[\begin{matrix} 0&0&0\\ 0&-1 &0\\ \nu_1&0 &0\end{matrix}\right] , \ \ C^T = \left[\begin{matrix} 1&0&0\\0&1&0\\ \nu_2&0&1\end{matrix}\right], \ \ Y^T = \left[\begin{matrix} 0&0&0\\ 1&0&0\\ 0 &0&0\end{matrix}\right]
$$
Here $\nu_1\ne 0$. In this case connected group $H$ is given by 
$$
H = \left\{\left[\begin{matrix} b& t &-\nu_1\ln a + \nu_2b \ln b \\ 0& ab&0\\ 0&0&b\end{matrix}\right] : a>0, b>0, t \in \R\right\}.
$$
An irreducibly admissible group is here given by  $F H$ where
$$
F = \left\{ \left[\begin{matrix} \epsilon & 0&0\\ 0& 1 &0\\ 0&0&1\end{matrix}\right] : \epsilon \in \{ -1, 1\}\right\}.
$$


 \medskip
 
{\rm (3)} If the dimension of the characteristic space $\E$ is two while the characteristic functions are linearly dependent, then   
 $$
A^T =\left[\begin{matrix} 1&0&0\\ 0&0&0\\ 0&0 &0 \end{matrix}\right] , \ \ C^T = \left[\begin{matrix} 1&0&0\\ 0&1&0\\ 0 &0&\l \end{matrix}\right], \ \ Y^T = \left[\begin{matrix} 0&0&0\\ 1&0&0\\ 0 &0&0\end{matrix}\right].
$$
In this case $H$ is given by 
$$
H = \left\{ \left[\begin{matrix} ab & t &0\\ 0& b&0\\ 0&0&b^\l\end{matrix}\right] : a>0, b>0, t \in \R\right\}~,
$$
with $\lambda \not= 0$. 
An irreducibly admissible group is here given by  $F H$ where
$$
F = \left\{ \left[\begin{matrix} \epsilon_1 & 0&0\\ 0& 1 &0\\ 0&0&\epsilon_2\end{matrix}\right] : \epsilon_1, \epsilon_2 \in \{ -1, 1\}\right\}.
$$


 

\medskip

{\rm (4)} If the dimension of the space $\E$ generated by the characteristic vectors is one, then 
$$
A^T =\left[\begin{matrix} 1&0&0\\ 0&0&0\\ 0&\delta_1 &0\end{matrix}\right], \ \ \ C^T = \left[\begin{matrix} 1 &0&0\\ 0&1&0\\ 0&\delta_2  &1 \end{matrix}\right], \ \ \ Y^T = \left[\begin{matrix}0&0&0\\ 0&0&0\\ 1 &0&0\end{matrix}\right],  
$$
where $\delta_1^2 +  \delta_2^2 \ne 0.$
  In this case 
$$
H = \left\{ \left[\begin{matrix} ab&0 &t\\ 0& b &\delta_1  \ln a  + \delta_2 b \ln b\\ 0&0&b\end{matrix}\right] : a>0, b>0, t \in \R\right\}.
$$
An irreducibly admissible group is here given by  $F H$ where
$$
F = \left\{ \left[\begin{matrix} \epsilon_1 & 0&0\\ 0& \epsilon_2 &0\\ 0&0&1 \end{matrix}\right] : \epsilon_1, \epsilon_2 \in \{ -1, 1\}\right\}.
$$

\end{proposition}

\section{The non-solvable case}

\label{sect:nonsolv}

It remains to show that the only non-solvable connected matrix group $H< {\rm GL}(3,\mathbb{R})$ acting with open orbits and associated compact stabilizers is conjugate to $\mathbb{R}^+ \times SO(3)$. 
Referring to the definition in \cite{knapp}, a real \textbf{linear connected reductive group} is a closed connected group of real matrices which is invariant under transpose. A \textbf{real linear connected semisimple group} is a connected linear real reductive group with finite center. It is well-known that if $G$ is connected and semisimple, then its algebra must be semisimple (is simple or is a direct sum of simple algebras) and if $G$ is a connected reductive Lie group then its Lie algebra can be written as a direct sum of the central ideal with the commutator ideal which is necessarily semisimple.  \vskip 0.2cm 

Let $G$ be a semisimple connected Lie group of real matrices with Lie algebra
$\mathfrak{g.}$ Let $\Theta$ be the inverse transpose which defines an
automorphism on $G.$ $\Theta$ is the so called Cartan involution. Next, let $$T=\left\{  g\in G:\Theta\left(  g\right)  =g\right\}.$$  Then $T$ is a maximal compact subgroup of $G.$ We define an involution $\iota:\mathfrak{g}\rightarrow\mathfrak{g}$ such that $\iota\left(  X\right)
=-X^{T}.$ This involution is the derivative of $\Theta$ at the identity of the group. It is also a Lie algebra automorphism. As a linear operator defined over $\mathfrak{g},$ the map $\iota$ has two
distinct eigenvalues:\ $1,-1.$ Let $\mathfrak{k}$ be the eigenspace
corresponding to the eigenvalue $1$ and let $\mathfrak{p}$ be the eigenspace
corresponding to $-1.$ Next, the Cartan decomposition is defined
as $\mathfrak{g}=\mathfrak{k+p}.$ It is worth noting that the algebra
$\mathfrak{k}$ contains skew-symmetric matrices, while $\mathfrak{p}$ consists
of symmetric matrices.
\vskip 0.2cm \noindent The following facts will be important.

\begin{enumerate}
\item The map $T\times\mathfrak{p}\rightarrow G$ given by $\left(  z,X\right)
\mapsto z\exp X$ is a diffeomorphism onto $G.$ For a proof, we refer to \cite[Page $362$]{knapp}.

\item Let $\mathfrak{h}$ be a four-dimensional Lie algebra which is not solvable. Let $\mathfrak{h}=\mathfrak{r}+\mathfrak{s}$ be the Levi decomposition of the algebra where $\mathfrak{s}$ is the semisimple part of the decomposition and $\mathfrak{r}$ is the solvable part. In this case, the dimension of the semisimple part must be equal to three. Next let $\mathfrak{k+p}$ be the Cartan decomposition of $\mathfrak{s}.$ If the dimension of the algebra $\mathfrak{k}$ is equal to one then the algebra $\mathfrak{s}$ is isomorphic to  $\mathfrak{sl}(2,\mathbb{R})$. Otherwise, $\mathfrak{s}$ is isomorphic to $\mathfrak{so}(3,\mathbb{R}).$
In fact it is well-known that any non-solvable four dimensional Lie algebra is Lie isomorphic to either $\mathfrak{sl}(2,\mathbb{R})\oplus\mathbb{R}$ or it is Lie isomorphic to $\mathfrak{so}(3,\mathbb{R})\oplus\mathbb{R}$; see the class labelled $U3$ in \cite{MacCallum}. 

\item Up to isomorphism there is one $(m+1)$-dimensional irreducible representation
of $\mathfrak{sl}(2,\mathbb{R})$ for each natural number $m,$ and every
finite-dimensional representation of $\mathfrak{sl}(2,\mathbb{R})$ is
completely reducible (thanks to Weyl's theorem.)
\end{enumerate}

Let $E_3$ be the identity matrix of order three. We shall now prove that if $H$ is a four-dimensional non-solvable Lie group of matrices of order three which acts with open orbits, then up to conjugation $H$ is one of the following groups
\begin{itemize}
\item $\mathrm{SO}_0(2,1)\cdot \exp(\mathbb{R}E_3)$
\item $\mathrm{SO}(3,\mathbb{R})\cdot \exp(\mathbb{R}E_3)$
\end{itemize}
The first group will then be excluded on the grounds that one of its open orbits has noncompact stabilizers. 

\noindent\textbf{Case (A)} Let us choose the following basis for $\mathfrak{sl}(2,\mathbb{R})$
satisfying
\[
\left[  H,X\right]  =2X,\left[  H,Y\right]  =2Y,\left[  X,Y\right]  =H.
\]
Let $\pi$ be an irreducible representation of this algebra in $\mathfrak{gl}\left(  V\right)  $ where $V$ is a three-dimensional vector space. Let
$\alpha$ be a real eigenvalue for the element $\pi\left(  H\right)  $ (there must be one since the dimension of $V$ is odd),
and let $v$ be a corresponding
eigenvector. By the irreducibility of the representation, the set $\left\{
v,\pi\left(  X\right)  v,\pi\left(  Y\right)  v\right\}  $ necessarily spans
$V.$ Moreover, using the fact that $\left[  H,X\right]  =2X,\left[  H,Y\right]  =2Y,$ it is easy to see that the endomorphism $\pi(H)$ is diagonalizable with eigenvalues $\alpha,\alpha+2,\alpha-2$. Next, let $w$ be
an eigenvector with corresponding eigenvalue $\alpha-2$ for $\pi(H).$ Using the structure
of the Lie algebra, with straightforward calculations, we check that the span
of the set $\left\{  w,\pi(X)w,\pi(X)^{2}w\right\}$ is invariant under the
action of $\pi(\mathfrak{h}).$ Appealing to the irreducibility of $\pi$ the
elements $w,\pi(X)w,\frac{1}{2}\pi(X)^{2}w$ form a basis for the vector space
$V.$ With respect to this basis, we obtain
\[
\pi\left(  H\right)  =\left[
\begin{array}
[c]{ccc}
\alpha-2 & 0 & 0\\
0 & \alpha & 0\\
0 & 0 & \alpha+2
\end{array}
\right]  ,\pi\left(  X\right)  =\left[
\begin{array}
[c]{ccc}
0 & 0 & 0\\
1 & 0 & 0\\
0 & 2 & 0
\end{array}
\right]  ,\pi\left(  Y\right)  =\left[
\begin{array}
[c]{ccc}
0 & \mu_{1} & 0\\
0 & 0 & \mu_{2}\\
0 & 0 & 0
\end{array}
\right]
\]
for some real numbers $\mu_{1},\mu_{2}.$ Next since $\left[
X,Y\right]  =H,$ it immediately follows that
\[
\pi\left(  H\right)  =\left[
\begin{array}
[c]{ccc}
-2 & 0 & 0\\
0 & 0 & 0\\
0 & 0 & 2
\end{array}
\right]  ,\pi\left(  X\right)  =\left[
\begin{array}
[c]{ccc}
0 & 0 & 0\\
1 & 0 & 0\\
0 & 2 & 0
\end{array}
\right]  ,\pi\left(  Y\right)  =\left[
\begin{array}
[c]{ccc}
0 & 2 & 0\\
0 & 0 & 1\\
0 & 0 & 0
\end{array}
\right]  .
\]
Furthermore, up to isomorphism the irreducible representation $\pi$ is unique. Let $\tau$ be a three-dimensional faithful representation of a Lie algebra
with basis $H,X,Y,Z$ satisfying the following non-trivial Lie brackets
\[
\left[  H,X\right]  =-2X,\left[  H,Y\right]  =2Y,\left[  Y,X\right]  =H
\]
such that the real span of $H,X,Y$ is $\mathfrak{sl}\left(  2,\mathbb{R}\right).$ This is a reductive Lie algebra with a one-dimensional center spanned by the element $Z$. Next, it is clear that if the representation $\tau$ is irreducible then up to isomorphism
\[
\tau\left(  H\right)  =\left[
\begin{array}
[c]{ccc}%
2 & 0 & 0\\
0 & 0 & 0\\
0 & 0 & -2
\end{array}
\right]  ,\tau\left(  X\right)  =\left[
\begin{array}
[c]{ccc}%
0 & 0 & 0\\
1 & 0 & 0\\
0 & 2 & 0
\end{array}
\right]  ,\tau\left(  Y\right)  =\left[
\begin{array}
[c]{ccc}%
0 & 2 & 0\\
0 & 0 & 1\\
0 & 0 & 0
\end{array}
\right]  ,\tau\left(  Z\right)  =\left[
\begin{array}
[c]{ccc}%
1 & 0 & 0\\
0 & 1 & 0\\
0 & 0 & 1
\end{array}
\right]  .
\]
Let $\mathfrak{h}$ be the linear span of the matrices above. We shall show
that this algebra is isomorphic to $\mathfrak{so}(2,1)\oplus\mathbb{R}$ where
\[
\mathfrak{so}(2,1)=\left\{  M\in\mathfrak{sl}(3,\mathbb{R}):M\left[
\begin{array}
[c]{ccc}%
1 & 0 & 0\\
0 & 1 & 0\\
0 & 0 & -1
\end{array}
\right]  +\left[
\begin{array}
[c]{ccc}%
1 & 0 & 0\\
0 & 1 & 0\\
0 & 0 & -1
\end{array}
\right]  M^{T}=0\right\}  .
\]
Put
\[
J=\left[
\begin{array}
[c]{ccc}%
0 & 1 & 0\\
-\frac{1}{2} & 0 & \frac{1}{2}\\
\frac{1}{2} & 0 & \frac{1}{2}%
\end{array}
\right]  =\left[
\begin{array}
[c]{ccc}%
0 & -1 & 1\\
1 & 0 & 0\\
0 & 1 & 1
\end{array}
\right]  ^{-1}.
\]
Indeed, with straightforward computations, we verify that
\[
J\left(  \tau\left(  -\frac{Y}{2}\right)  +\tau\left(  \frac{X}{2}\right)
\right)  J^{-1}=\left[
\begin{array}
[c]{ccc}%
0 & -1 & 0\\
1 & 0 & 0\\
0 & 0 & 0
\end{array}
\right]  ,J\tau\left(  \frac{H}{2}\right)  J^{-1}=\left[
\begin{array}
[c]{ccc}%
0 & 0 & 0\\
0 & 0 & -1\\
0 & -1 & 0
\end{array}
\right]
\]
and
\[
J\left(  \tau\left(  \frac{Y}{2}\right)  +\tau\left(  \frac{X}{2}\right)
\right)  J^{-1}=\left[
\begin{array}
[c]{ccc}%
0 & 0 & 1\\
0 & 0 & 0\\
1 & 0 & 0
\end{array}
\right]  .
\]
We now define $\mathfrak{h}_{1}$ to be the algebra spanned by $\tau
(H),\tau(X),\tau(Y).$ It follows from the computations above that
$\mathfrak{h}_{1}$ is isomorphic to $\mathfrak{so}(2,1)$ which is spanned by
\[
K=\left[
\begin{array}
[c]{ccc}%
0 & -1 & 0\\
1 & 0 & 0\\
0 & 0 & 0
\end{array}
\right]  ,A=\left[
\begin{array}
[c]{ccc}%
0 & 0 & 1\\
0 & 0 & 0\\
1 & 0 & 0
\end{array}
\right]  ,B=\left[
\begin{array}
[c]{ccc}%
0 & 0 & 0\\
0 & 0 & -1\\
0 & -1 & 0
\end{array}
\right]  .
\]
Indeed it is easy to check that for any real numbers $t,r,s$
\[
\left(  tK+rA+sB\right)  \left[
\begin{array}
[c]{ccc}%
1 & 0 & 0\\
0 & 1 & 0\\
0 & 0 & -1
\end{array}
\right]  +\left[
\begin{array}
[c]{ccc}%
1 & 0 & 0\\
0 & 1 & 0\\
0 & 0 & -1
\end{array}
\right]  \left(  tK+rA+sB\right)  ^{T}
\]
is equal to the zero matrix. Thus, $\mathfrak{h}$ is isomorphic to
$\mathfrak{so}(2,1)\oplus\mathbb{R}$.  Next, we aim to compute the dual action
of $SO_{0}(2,1)\exp\mathbb{R}E_{3}$ (its Lie algebra is spanned by
$K,A,B,E_{3}$) where $E_{3}$ is the identity matrix of order three. The
corresponding Cartan decomposition is given by $\mathfrak{so}\left(
2,1\right)  =\mathfrak{k+p}$ where $\mathfrak{k}=\mathbb{R}K$ and
$\mathfrak{p}=\mathbb{R}A+\mathbb{R}B.$ Let $\mathbb{T}$ be a maximal torus
subgroup with Lie algebra $\mathfrak{k.}$ Since the map $\mathbb{T}\times\mathfrak{p}\rightarrow SO_{0}\left(  2,1\right)  $ where $\left(
z,X\right)  \mapsto z\exp X$ is a diffeomorphism then every element of
$SO_{0}\left(  2,1\right)  $ can be written as
\[
\left[
\begin{array}
[c]{ccc}%
\cos\theta & -\sin\theta & 0\\
\sin\theta & \cos\theta & 0\\
0 & 0 & 1
\end{array}
\right]  \exp\left(  aA+bB\right)
\]
where $\left(  \theta,a,b\right)  \in\left[  0,2\pi\right)  \times
\mathbb{R}\times\mathbb{R}$. Writing $\left(  a,b\right)  $ in polar
coordinates such that $a=r\cos s,$ and $b=r\sin s,$ we obtain that
\[
\exp\left(  aA+bB\right)  =\exp\left[
\begin{array}
[c]{ccc}%
0 & 0 & r\cos s\\
0 & 0 & r\sin s\\
r\cos s & r\sin s & 0
\end{array}
\right]
\]
is equal to
\[
\left[
\begin{array}
[c]{ccc}%
\cos^{2}s\cosh r+(\sin s)^{2} & -\sin2s\sinh^{2}(\frac{r}{2}) & \cos s\sinh
r\\
-\sin2s\sinh^{2}(\frac{r}{2}) & \cos^{2}s+\cosh r(\sin s)^{2} & -\sin s\sinh
r\\
\cos s\sinh r & -\sin s\sinh r & \cosh r
\end{array}
\right]  .
\]
Thus, every element of $SO_{0}(2,1)\exp\mathbb{R}E_{3}$ can be uniquely
written as
\[
\left[
\begin{array}
[c]{ccc}%
\cos\theta & -\sin\theta & 0\\
\sin\theta & \cos\theta & 0\\
0 & 0 & 1
\end{array}
\right]  \left[
\begin{array}
[c]{ccc}%
\cos^{2}s\cosh r+(\sin s)^{2} & -\sin2s\sinh^{2}(\frac{r}{2}) & \cos s\sinh
r\\
-\sin2s\sinh^{2}(\frac{r}{2}) & \cos^{2}s+\cosh r\sin^{2}s & -\sin s\sinh r\\
\cos s\sinh r & -\sin s\sinh r & \cosh r
\end{array}
\right]  \left[
\begin{array}
[c]{ccc}%
e^{t} & 0 & 0\\
0 & e^{t} & 0\\
0 & 0 & e^{t}%
\end{array}
\right]
\]
and the dual action of the group can now be computed in a very precise manner.
$\allowbreak$Let
\[
\Omega^{\star}=\left\{  \left[
\begin{array}
[c]{c}%
x\\
y\\
z
\end{array}
\right]  :x^{2}+y^{2}-z^{2}\neq0\right\}
\]
be a Zariski open subset of the dual of $V.$ We shall prove that
$\Omega^{\star}$ is $SO_{0}\left(  2,1\right)  \exp\mathbb{R}E_{3}$-invariant
and
\[
\Sigma=\left\{  \omega,\xi^{+},\xi^{-}\right\}  \text{ where }\omega=\left[
\begin{array}
[c]{c}%
0\\
1\\
0
\end{array}
\right]  ,\xi^{+}=\left[
\begin{array}
[c]{c}%
0\\
0\\
1
\end{array}
\right]  ,\xi^{-}=\left[
\begin{array}
[c]{c}%
0\\
0\\
-1
\end{array}
\right]
\]
is a cross-section for the dual orbits in $\Omega^{\star}$. First, we check
that the dual orbit of $\xi^{\pm}$ is equal to
\[
\Omega^{\pm}=\left\{  \left[
\begin{array}
[c]{c}%
e^{t}\sin s\sinh r\\
-e^{t}\cos s\sinh r\\
\pm e^{t}\cosh r
\end{array}
\right]  :t,r,s\in\mathbb{R}\right\}  .
\]
Next, let $\Omega=\Omega^{+}\cup\Omega^{-}.$ Then
\[
\Omega={\bigcup\limits_{h>0}}\left\{  \left[
\begin{array}
[c]{c}%
hx\\
hy\\
hz
\end{array}
\right]  :x^{2}+y^{2}-z^{2}=-1\right\}  =\left\{  \left[
\begin{array}
[c]{c}%
x\\
y\\
z
\end{array}
\right]  :x^{2}+y^{2}-z^{2}<0\right\}  .
\]
$\allowbreak$ The corresponding stabilizer subgroup for $\xi^{\pm}$ in
$SO_{0}(2,1)\exp\mathbb{R}E_{3}$ is equal to a maximal torus subgroup of
$SO_{0}(2,1).$ Thus the stabilizer subgroup for every element in $\Omega$ is
compact. Finally, the orbit of $\omega$ is given by $\allowbreak$%
\[
\Omega^{\prime}=\left\{  \left[
\begin{array}
[c]{c}%
x\\
y\\
z
\end{array}
\right]  :x^{2}+y^{2}-z^{2}>0\right\}  .
\]
However, every element of the one-parameter group $\exp\left(  \mathbb{R}%
A\right)  $ stabilizes $\omega.$ Finally, our claim follows from the fact that
$\Omega^{\star}=\Omega\cup\Omega^{\prime}.$ Next, let $\pi$ be a
two-dimensional irreducible representation of $\mathfrak{sl}(2,\mathbb{R)}$
acting on a vector space $V.$ Since the algebra generated by $\pi(H),\pi(X)$
is solvable then there exists a basis for the vector space $V$ such that
\[
\pi\left(  H\right)  =\left[
\begin{array}
[c]{cc}%
\alpha & 0\\
0 & \alpha+2
\end{array}
\right]  ,\pi\left(  X\right)  =\left[
\begin{array}
[c]{cc}%
0 & 0\\
1 & 0
\end{array}
\right]
\]
and $v$ is an eigenvector for $\pi\left(  H\right)  $ with corresponding
eigenvalue $\alpha.$ Now, since $\left\{  v,\pi\left(  X\right)  v\right\}  $
is a linearly independent set, it is a basis. With respect to this basis, we
obtain
\[
\pi\left(  H\right)  =\left[
\begin{array}
[c]{cc}%
-1 & 0\\
0 & 1
\end{array}
\right]  ,\pi\left(  X\right)  =\left[
\begin{array}
[c]{cc}%
0 & 0\\
1 & 0
\end{array}
\right]  ,\pi\left(  Y\right)  =\left[
\begin{array}
[c]{cc}%
0 & 1\\
0 & 0
\end{array}
\right]
\]
and clearly%
\[
\pi\left(  Y\right)  -\pi\left(  X\right)  =\left[
\begin{array}
[c]{cc}%
0 & 1\\
-1 & 0
\end{array}
\right]  ,\left[
\begin{array}
[c]{cc}%
1 & 0\\
0 & -1
\end{array}
\right]  ,\left[
\begin{array}
[c]{cc}%
0 & 1\\
0 & 0
\end{array}
\right]
\]
form a basis for $\pi\left(  \mathfrak{sl}(2,\mathbb{R)}\right)  .$ Let $\tau$
be a three-dimensional reducible representation of $\mathfrak{sl}\left(
2,\mathbb{R}\right)  \oplus\mathbb{R}$. Then up to conjugation, the image of
$\mathfrak{sl}\left(  2,\mathbb{R}\right)  \oplus\mathbb{R}$ via this
representation is spanned by matrices
\[
K=\left[
\begin{array}
[c]{ccc}%
0 & 1 & 0\\
-1 & 0 & 0\\
0 & 0 & 0
\end{array}
\right]  ,A=\left[
\begin{array}
[c]{ccc}%
1 & 0 & 0\\
0 & -1 & 0\\
0 & 0 & 0
\end{array}
\right]  ,N=\left[
\begin{array}
[c]{ccc}%
0 & 1 & 0\\
0 & 0 & 0\\
0 & 0 & 0
\end{array}
\right]  ,D=\left[
\begin{array}
[c]{ccc}%
\lambda & 0 & 0\\
0 & \lambda & 0\\
0 & 0 & \beta
\end{array}
\right]
\]
where $\lambda,\beta$ are some real numbers. We observe that the matrices
$K,A,N$ span an algebra which is isomorphic to $\mathfrak{sl}\left(
2,\mathbb{R}\right)  $. Next, we consider the dual action of $H=SL\left(
2,\mathbb{R}\right)  \exp\left(  \mathbb{R}D\right)  .$ Using the Cartan
decomposition of $\mathfrak{sl}\left(  2,\mathbb{R}\right)  ,$ we write
elements of $H$ as $\exp\left(  \theta K\right)  \exp\left(  t_{1}%
A+t_{2}N+t_{3}D\right)  $ where $\theta\in\lbrack0,2\pi),$ and $t_{1}%
,t_{2},t_{3}$ are some real numbers. Assuming that $\beta\neq0,$ there are two
open orbits, and a cross-section for these orbits is given by
\[
\Sigma=\left\{  \xi^{+},\xi^{-}\right\}  \text{ where }\xi^{+}=\left[
\begin{array}
[c]{c}%
1\\
0\\
1
\end{array}
\right]  ,\xi^{-}=\left[
\begin{array}
[c]{c}%
1\\
0\\
-1
\end{array}
\right]  .
\]
The stabilizer subgroups corresponding to the elements $\xi^{+},\xi^{-}$ are
not trivial and not contained in the maximal torus subgroup of $H$ (and thus
not compact.) Indeed the one-parameter group
\[
\left\{  \left[
\begin{array}
[c]{ccc}%
1 & 0 & 0\\
t & 1 & 0\\
0 & 0 & 1
\end{array}
\right]  :t\in\mathbb{R}\right\}
\]
stabilizes $\xi^{+}$ and $\xi^{-}.$ Thus the dilation group $SL\left(
2,\mathbb{R}\right)  \exp\left(  \mathbb{R}D\right)  $ is not irreducibly
admissible. Assume that $\beta=0$ (thus $\lambda\neq0.$) There are no open
orbits associated with this dilation group. Thus the corresponding dilation
groups are not irreducibly admissible. \vskip 0.5cm

\noindent\textbf{Case (B)} Let $\mathfrak{h}$ be a Lie algebra of matrices of order three which is
isomorphic to the direct sum of $\mathfrak{so}\left(  2,\mathbb{R}\right)
\oplus\mathbb{R}$. Up to conjugation, $\mathfrak{h}$ is spanned by matrices
\[
X_{1}=\left[
\begin{array}
[c]{ccc}%
0 & -1 & 0\\
1 & 0 & 0\\
0 & 0 & 0
\end{array}
\right]  ,X_{2}=\left[
\begin{array}
[c]{ccc}%
0 & 0 & 0\\
0 & 0 & -1\\
0 & 1 & 0
\end{array}
\right]  ,X_{3}=\left[
\begin{array}
[c]{ccc}%
0 & 0 & -1\\
0 & 0 & 0\\
1 & 0 & 0
\end{array}
\right]  ,X_{4}=\left[
\begin{array}
[c]{ccc}%
1 & 0 & 0\\
0 & 1 & 0\\
0 & 0 & 1
\end{array}
\right]
\]
Now, let $H$ be a Lie group with Lie algebra $\mathfrak{h}$ and let
$H_{0}$ be the connected component of the identity. Then $H_{0}=SO\left(
3,\mathbb{R}\right)  \exp\left(  \mathbb{R}X_{4}\right)  .$ Next, let $v$ be a
unit vector. Clearly, there exists a pair of angles $\left(  \theta
,\phi\right)  $ such that $v=\exp\left(  \theta X_{1}\right)  \exp\left(  \phi
X_{2}\right)  e_{3}.$ Thus, for an arbitrary matrix $M$ in the orthogonal
group $SO\left(  3,\mathbb{R}\right)  $ such that $Me_{3}=v,$ $Me_{3}=\exp\left(  \theta X_{1}\right)  \exp\left(  \phi X_{2}\right)  e_{3}$ and
there exists an angle $\omega$ such that
\[
M=\exp\left(  \theta X_{1}\right)  \exp\left(  \phi X_{2}\right)  \exp\left(
\omega X_{1}\right)  .
\]
Thus an element of $SO\left(  3,\mathbb{R}\right)  \exp\left(  \mathbb{R}X_{4}\right)  $ can be written as follows:
\[
\left[
\begin{array}
[c]{ccc}%
\cos\theta & -\sin\theta & 0\\
\sin\theta & \cos\theta & 0\\
0 & 0 & 1
\end{array}
\right]  \left[
\begin{array}
[c]{ccc}%
1 & 0 & 0\\
0 & \cos\phi & -\sin\phi\\
0 & \sin\phi & \cos\phi
\end{array}
\right]  \left[
\begin{array}
[c]{ccc}%
\cos\omega & -\sin\omega & 0\\
\sin\omega & \cos\omega & 0\\
0 & 0 & 1
\end{array}
\right]  \left[
\begin{array}
[c]{ccc}%
e^{t} & 0 & 0\\
0 & e^{t} & 0\\
0 & 0 & e^{t}%
\end{array}
\right]  .
\]
Computing the dual action of $SO\left(  3,\mathbb{R}\right)  \exp\left(
\mathbb{R}X_{4}\right)  $, we obtain a single open orbit (the punctured space)
and the stabilizer for every element in $V^{\ast}$ is compact. Thus, this dilation group is irreducibly admissible.

\section{Compactly supported wavelets and atoms}

\label{sect:atom}

In this section we wish to establish the existence of compactly supported admissible wavelets for all irreducibly admissible subgroups  $H<{\rm GL}(3,\mathbb{R})$. For most groups, we will be able to show a much stronger statement, namely the existence of compactly supported {\em atoms} $\psi$ associated to $H$. These functions are distinguished by the fact that for all sufficiently dense lattices $\Gamma \subset \mathbb{R}^d$ and all sufficiently uniformly dense and uniformly discrete subsets $\Lambda \subset H$, the family $(\pi(hx,h) \psi)_{x \in \Gamma, h \in \Lambda}$ is a frame, and furthermore the decay of the frame coefficients of $f \in {\rm L}^2(\mathbb{R}^d)$ characterizes whether $f$ belongs to a suitable {\bf coorbit space} $Co(Y)$, where $Y$ is a Banach space of functions on $G$. 

It is known that suitably chosen bandlimited Schwartz functions can serve as atoms \cite{Fu_coorbit}, however the existence of atoms that are compactly supported {\em in space} is not generally established. In fact atoms are admissible vectors, i.e., $\mathcal{W}_\psi$ is a well-defined isometry for each admissible vector $\psi$, and it is currently not even known whether there exist compactly supported vectors for every irreducibly admissible matrix group. The recent papers \cite{Fu_atom,Fu_RT} contain sufficient criteria for the existence of such atoms, and below we will prove that these criteria are fulfilled {\em for all cases except case (2.l)}. Note that this does not exclude the existence of compactly supported atoms for this case, just that the methods from \cite{Fu_atom, Fu_RT} do not suffice for a positive answer. For the exceptional case (2.l), we will however establish the existence of a compactly supported admissible vector.

\subsection{Coorbit spaces and atomic decompositions}

\label{subsect:atom}

Coorbit space theory addresses two important questions in generalized wavelet analysis:  The consistent definition of spaces of sparse signals, and discretization of the continuous wavelet transforms. 
By definition, sparse signals should be described by fast decay of their wavelet transforms. One possible way of quantifying this decay behavior is to impose a weighted mixed ${\rm L}^p$-norm on the wavelet coefficients: Nice signals are those for which the wavelet coefficient decay is sufficiently fast to guarantee weighted integrability. One of the starting points of coorbit space theory was the realization that the scale of Besov spaces can be understood in precisely these terms \cite{FeiGr1}.

Coorbit space theory provides the same type of results for a much larger setup, including continuous wavelet transforms associated to irreducibly admissible matrix groups $H$. We will now briefly describe the main ingredients of this theory, for a general irreducibly admissible matrix group $H < {\rm GL}(d,\mathbb{R})$. We start out by fixing a Banach function space $Y$ on $G = \mathbb{R}^d \rtimes H$. A relevant choice from the approximation-theoretic point of view is ${\rm L}^p(G)$, with $1 \le p < 2$; see Remark \ref{rem:co_lp}. Furthermore, we assume that $Y$ is strongly invariant under both left and right translations of $G$, and that it is {\bf solid}: For every pair $F,H$ of measurable functions on $G$  with $F \in Y$ and $|H|\le |F|$ pointwise a.e., it follows that $H \in Y$ and $\|H \|_Y \le \|F \|$.  

 \begin{definition} \label{defn:weight}
 \begin{enumerate}
\item[(a)] A {\bf weight} on $G$ is a continuous function $v: G \to \mathbb{R}^+$ satisfying the submultiplicativity condition 
 \[
  v((x_1,h_1) (x_2,h_2) ) \le v(x_1,h_1) v(x_2,h_2)~.
 \]
\item[(b)]
A weight $v$ on $G$ is called {\bf polynomially bounded} if it satisfies
\[
 v(x,h) \le C (1+|x|)^s (1+\| h \|)^r~,
\] for suitable $r,s\ge0$.  
 \item[(c)]  Let $Y$ be a solid Banach function space on $G$.
   A weight $v_0$ is called {\bf control weight} for $Y$ if it is bounded from below, and satisfies 
\[
 v_0(x,h) = \Delta_G(x,h)^{-1} v_0((x,h)^{-1})~, 
\]
as well as 
 \[
\max \left( \|L_{(x,h)^{\pm 1}} \|_{Y \to Y},\| R_{(x,h)} \|_{Y \to Y},\|
R_{(x,h)^{-1}} \|_{Y \to Y} \Delta_G(x,h)^{-1} \right) \le v_0(x,h)
\] where $L_{(x,h)}$ and $R_{(x,h)}$ denote left and right translation by $(x,h) \in G$.
\end{enumerate}
 \end{definition}
 
 Note that two different notions of weights exist in this paper. Each is established terminology within its context, and for the remainder of the paper, only the definition of weights according to \ref{defn:weight} will play a role. For these reasons, we have refrained to further differentiate between the two, and prefer to live with a slight conflict of terminology as a consequence. 
 
 We next define the set of atoms. The role of the weight $v_0$ in the definition can be read as a {\rm compatibility condition} incurred by the choice of the function space $Y$: The faster $v_0$ grows, the more restrictive the definition of $\mathcal{B}_{v_0}$ becomes. This is comparable to the vanishing moment and smoothness conditions imposed on wavelet ONB's characterizing homogeneous Besov spaces $\dot{B}_{p,q}^\alpha$, which also depend on the parameters $p,q,\alpha$.  
 \begin{definition} \label{defn:atom}
  Let $Y$ denote a solid Banach function space on the locally compact group $G$, $U \subset G$ a compact neighborhood of the identity, and $F: G \to \mathbb{C}$. We let
  \[
   \left(\mathcal{M}_U^R F \right) (x) = \sup_{y \in U}|F(yx^{-1})|
  \] denote the right local maximum function of $F$ with respect to $U$. Given a weight $v_0$ on $G$, 
  let 
  \[ {\rm L}^1_{v_0}(G) = \{ F : G \to \mathbb{C}~: F v_0 \in {\rm L}^1(G) \} ~.\]
 Now $0 \not= \psi \in {\rm L}^2(\mathbb{R}^d)$ is called a {\bf $v_0$-atom}, if
  \[
\mathcal{M}_U^R (W_\psi \psi)  \in {\rm L}^1_{v_0}(G) ~,
  \]
  and we denote the set of all such $\psi$ by $\mathcal{B}_{v_0}$. 
   \end{definition}
 
 \begin{remark} \label{rem:atom_adm}
  We note that 
  \[
   |\mathcal{W}_\psi \psi| \le \mathcal{M}_U^R \left( \mathcal{W}_\psi \psi \right) 
  \] implies that for every $\psi \in \mathcal{B}_{v_0}$, 
  \[ \mathcal{W}_\psi \psi \in {\rm L}^1_{v_0}(G) \cap {\rm L}^\infty(G) \subset {\rm L}^2(G) ~,\]
  since $v_0$ is assumed bounded from below. 
  Hence every $\psi \in \mathcal{B}_{v_0}$ is an admissible vector. 
 \end{remark}

Now the central results of \cite{FeiGr0,FeiGr1,FeiGr2,Gr} with respect to consistency and discretization can be briefly summarized as follows:
\begin{enumerate}
\item Given a suitable Banach function space $Y$ on $G$ with control weight $v_0$, fix a nonzero $\psi \in \mathcal{A}_{v_0} := \{ \psi\in L^2(\R^d) : \mathcal W_\psi\psi \in L^1_{v_0}\}$, and define 
$\mathcal{H}_{1,v_0} = \{ f \in {\rm L}^2(\mathbb{R}^d): \mathcal{W}_\psi f \in  {\rm L}^1_{v_0} \}$. 
Let $\mathcal{H}_{1,v_0}^{\sim}$ denote the space of conjugate-linear bounded functionals on $\mathcal{H}_{1,v_0}$. Then the wavelet transform $\mathcal{W}_\psi f$ of $f \in \mathcal{H}_{1,v_0}^{\sim}$ can be defined by canonical extension of the formula for ${\rm L}^2$-functions.  
\item Fix a nonzero $\psi \in \mathcal{B}_{v_0}$, and define 
\[
 {\rm Co} Y = \{ f \in \mathcal{H}_{1,v_0}^{\sim}~:~ \mathcal{W}_\psi f \in Y \}~, 
\] with the obvious norm $\| f \|_{{\rm Co} Y} = \| \mathcal{W}_\psi f \|_Y$. Then ${\rm Co} Y$ is a well-defined Banach space, and it is {\em independent} of the choice of $\psi \in \mathcal{B}_{v_0} \setminus \{ 0 \}$: Changing $\psi$ results in an equivalent norm, and in the same coorbit space $Co Y$.
\item Fix a nonzero $\psi \in \mathcal{B}_{v_0}$. Then the coorbit space norm is equivalent to the discretized norm $\left\| \mathcal{W}_\psi f|_Z \right\|_{Y_d}$, for all suitably dense and discrete subsets $Z \subset G$, with a suitably defined Banach sequence space $Y_d$. This also gives rise to atomic decompositions, i.e., discrete systems of wavelets that provide frame decompositions converging not just in ${\rm L}^2$, but also in the coorbit space norms. 
\end{enumerate}
 
\begin{remark} \label{rem:co_lp}
Readers who prefer explicitness over generality may take the spaces $Y = {\rm L}^p(G)$ as template for the general setup, with $p \in [1,\infty]$. Here the above statements can be further motivated by non-linear approximation. To begin with, the associated control weight can be taken as $v_0(x,h) = \max(1,\Delta_G(h))$, valid for all $p \in [1,\infty]$ \cite[Lemma 2.3 and subsequent remarks]{Fu_atom}. The benefit of using a common weight is that it gives rise to atomic decompositions that converge simultaneously in a whole range of coorbit spaces.

We will now shortly explain the relevance of $Co({\rm L}^p(G))$ from the point of view of non-linear approximation: Let $1 \le p < 2$, and let $\psi \in \mathcal{B}_{v_0}$. Then there exists a discrete subset $\Lambda \subset G$ such that  $(\pi(\lambda) \psi)_{\lambda \in \Lambda}$ is a frame of ${\rm L}^2(G)$, which means that there exist constants $0 < A \le B$, depending on $
\Lambda$ and $\psi$, such that for all $f \in {\rm L}^2(\mathbb{R}^3)$, 
\[
 A \| f\|_2^2 \le \sum_{(x,h) \in 
 \Lambda} |\langle f, \pi(hx,h) \psi \rangle|^2 \le B \| f \|_2^2~.
\] In this case, there exists a dual frame that can be used to write down a discrete version of the wavelet inversion formula, converging unconditionally in ${\rm L}^2$.

But atoms allow much more than just frame expansions converging in ${\rm L}^2$. As a matter of fact, the elements $f \in Co({\rm L}^p)$ are characterized by any of the following equivalent statements:
\begin{itemize}
 \item[(a)] The discrete coefficient sequence $(\langle f, \pi(\lambda) \psi)_{\lambda \in \Lambda}$ is $p$-summable. 
 \item[(b)] $f$ has an expansion
 \[
  f = \sum_{(x,h) \in \Lambda} c_{\lambda} \pi(\lambda) \psi ~,
 \] with $p$-summable coefficients $(c_\lambda)_{\lambda\in \Lambda}$. 
\end{itemize}
The second property has interesting consequences for nonlinear approximation: If we let
\begin{equation} \label{eqn:approx_error}
 E_n(f; (\pi(\lambda) \psi)_{\lambda}) = \inf \left\{ \left\| f - \sum_{\lambda \in \Lambda'} c_\lambda \pi(\lambda) \psi  \right\|_{\mathcal{H}} ~:~ c_\lambda \in \mathbb{C}~,~|\Lambda'| \le n \right\} ~.\end{equation}
 then approximating $f$ using the $n$ largest coefficients in (b) leads to the following estimate for the approximation error
 \[
  E_n(f, (\pi(\lambda) \psi)_{\lambda}) \le C_\epsilon \| f \|_{Co({\rm L}^p)} n^{-p/2+\epsilon}~,
 \] valid for any $\epsilon >0$ and $f \in Co({\rm L}^p)$. We refer to \cite[Section 1.1]{Fu_atom} for more details.
These observations clearly justify the interpretation of $Co({\rm L}^p)$ as a space of sparse signals. 
 
 It should also be mentioned that the coorbit scheme extends to the quasi-Banach setting; see \cite{Rau,VoDiss}. In particular, $Y = {\rm L}^p(G)$ with $p<1$ can also be employed, which leads to even sharper decay estimates for nonlinear approximation. Note however that the vanishing moment conditions for atoms that we will employ below have not yet been established for quasi-Banach spaces. 
\end{remark}

 Note that the discussion in this subsection is devoid of meaning if the set $\mathcal{B}_{v_0}$ is empty. Clearly, checking the condition in Definition \ref{defn:atom} poses a significant challenge, since it requires estimating Wiener amalgam space norms of continuous wavelet transforms.  However, for the setting of generalized wavelet transforms, the following result is obtained by combining \cite[Theorem 2.1]{Fu_coorbit} with \cite[Lemma 2.7]{Fu_coorbit}.
 \begin{theorem}
  Assume that $0 \not= \widehat{\psi} \in C_c^\infty(\mathcal{O})$. Then $\psi \in \mathcal{B}_{v_0}$. 
 \end{theorem}

\subsection{Vanishing moment conditions for admissible vectors and atoms}

\label{subsect:vm}

The previous subsection was mainly intended to motivate the interest in the set $\mathcal{B}_{v_0}$ of atoms. At this point the only easily accessible atoms are compactly supported in frequency, and the existence of atoms that are compactly supported in space is open. Sufficient criteria for this latter class have been formulated in \cite{Fu_atom,Fu_RT}, and we will now briefly recount these conditions. 

We denote the unique open dual orbit of $H$ by $\mathcal{O} \subset \mathbb{R}^d$, and its complement by $\mathcal{O}^c$.
We aim for statements that a function $\psi$ is admissible (or an atom) provided it has sufficient decay, sufficient smoothness, and sufficiently many vanishing moments in $\mathcal{O}^c$. Only the last condition actually reflects our specific choice of the matrix group $H$, and it is rigorously defined as follows: 
\begin{definition} \label{defn:van_mom}
 Let $r \in \mathbb{N}$ be given.
 $f \in {\rm L}^1(\mathbb{R}^d)$ {\bf has vanishing moments in $\mathcal{O}^c$ of order $r$}  if all
 distributional derivatives $\partial^\alpha \widehat{f}$ with $|\alpha|\le r$ are
 continuous functions, and all derivatives of degree $|\alpha|<r$ are vanishing on $\mathcal{O}^c$.
\end{definition} 
Note that under suitable integrability conditions on $\psi$, the vanishing moment conditions are equivalent to 
\[
 \forall |\alpha| < r,\forall \xi \in \mathcal{O}^c ~:~\int_{\mathbb{R}^d} x^\alpha \psi(x) e^{-2 \pi i \langle \xi, x \rangle} dx = 0~. 
\]
It is important to  note that \cite[Lemma 4.1]{Fu_coorbit} guarantees, for any $r>0$, the existence of $\psi \in C_c^\infty(\mathbb{R}^d)$ with vanishing moments of order $r$.  Hence sufficient vanishing moment criteria for atoms guarantee the existence of compactly supported atoms. 

The central tool for the formulation of vanishing moment conditions for admissible vectors and atoms is an auxiliary function $A: \mathcal{O} \to \mathbb{R}^+$ defined as follows:
Given any point $\xi \in \mathcal{O}$, let ${\rm
dist}(\xi,\mathcal{O}^c)$ denote the minimal distance of $\xi$ to
$\mathcal{O}^c$, and define
\[
 A(\xi) := \min \left( \frac{{\rm dist}(\xi,\mathcal{O}^c)}{1+\sqrt{|\xi|^2-{\rm dist}(\xi,\mathcal{O}^c)^2}},
 \frac{1}{1+|\xi|} \right)~. 
\]
By definition, $A$ is a continuous function with $A(\cdot) \le 1$, and it is extended continuously to all of $\mathbb{R}^d$ by setting it to zero on $\mathcal{O}^c$.
If $\eta \in \mathcal{O}^c$ denotes an element of minimal distance to $\xi$, the fact that $\mathbb{R}^+ \cdot \eta  \subset \mathcal{O}^c$ then entails that $\eta$ and $\xi-\eta$ are orthogonal with respect to the standard scalar product on $\mathbb{R}^d$, and we obtain the more transparent expression
\[
  A(\xi) = \min \left( \frac{|\xi-\eta|}{1+|\eta|},
 \frac{1}{1+|\xi|} \right)~. 
\]

We can now formulate sufficient criteria for the existence of compactly supported admissible functions and atoms. Note that
part (a) provides access to compactly supported atoms for coorbit spaces associated to weighted-mixed ${\rm L}^p$ spaces with 
polynomial weights. 
\begin{proposition} \label{prop:van_moments}
Assume that $v_0$ is a polynomially bounded weight on $G$. 
 Define $A_H : H \to \mathbb{R}^+$ by 
 \[
  A_H(h) = A(h^T \xi)~. 
 \]
 \begin{enumerate}
  \item[(a)]
   If there exists an exponent $e>1$ and a constant $C>0$ such that, for all $h \in H$,
   \[
    \Delta_G(h) A_H(h)^e \le C
   \] then $\psi \in C_c^\infty(\mathbb{R}^d)$ with vanishing moments of order $r > (d+e)/2$ is admissible.
   \item[(b)] 
 If there exists an exponent $e>1$ and a constant $C>0$ such that, for all $h \in H$,
   \[
    \| h^{\pm 1} \| A_H(h)^e \le C
   \] then there exists an $r>0$, explicitly computable from $v_0$ and $d$, such that every compactly supported smooth function sufficiently many vanishing moments of order $r$ is contained in $\mathcal{B}_{v_0}$.
 \end{enumerate}
\end{proposition}

Note that the condition in Part (b) does not depend on the precise choice of matrix norms. 

\begin{proof}

For part (a), we recall the general admissibility condition derived in \cite{Fu96}: Fixing a vector $\xi_0 \in \mathcal{O}$, the fact that the continuous homomorphism $\Delta_G$ restricted to the compact stabilizer $H_{\xi_0}$ is constant implies that letting $\Phi: \mathcal{O} \to \mathbb{R}^+$, with $\Phi(h^T \xi_0) =  \Delta_G(h)$,  yields a well-defined continuous function.
Furthermore, $\Phi$ has the property that $0 \not= \psi$ is admissible iff 
\[
 \int_{\mathcal{O}} |\widehat{\psi}(\xi)|^2 \Phi(\xi) d\xi < \infty~. 
\] 
Our assumption on $e$ and $C$ yield
\[
 \Phi(\xi) A(\xi)^e \le C~.
\]

If $\psi$ is a compactly supported $C^\infty$ functions with vanishing moments of order $r$, we get by \cite[Lemma 3.4]{Fu_coorbit}
that 
\[
 |\widehat{\psi}(\xi)| \le C' |\widehat{\psi}|_{r,m} A(\xi)^r~.
\] Here we used the Schwartz norm 
\[
| f |_{r,m} = \sup_{x \in \mathbb{R}^d, |\alpha| \le r} (1+|x|)^{m}
|\partial^\alpha f (x)|~.
\]  Hence we obtain
\begin{eqnarray*}
 \int_{\mathcal{O}} |\psi(\xi)|^2 \Phi(\xi) dx & \le & C' |\psi|_{r,m} \int_{\mathcal{O}} A(\xi)^{2r} \Phi(\xi) d\xi \\
 & \le & C C' |\psi|_{r,m}  \int_{\mathcal{O}} A(\xi)^{2r-e} d\xi \\
 & \le & C C' |\psi|_{r,m} \int_{\mathbb{R}^d} (1+|\xi|)^{e-2r} d\xi < \infty~,
\end{eqnarray*}
as soon as $2r -e > d$. 

Part (b) follows combining \cite[Theorem 3.4] {Fu_atom} with \cite[Lemma 2.8, Theorem 2.12]{Fu_RT}.
\end{proof}

\subsection{Proof of Theorem \ref{thm:main} (b), (c)}

We will verify the conditions of Proposition \ref{prop:van_moments}. First note that the existence of compactly supported admissible vectors or atoms is clearly conjugation invariant: If the admissible matrix groups $H_1$ and $H_2$ are related via $H_2 = g H_1 g^{-1}$, with $g \in {\rm GL}(3,\mathbb{R})$, then every such vector for $H_2$ is obtained from a vector with analogous properties by dilation with $g$, which preserves compact support.

Furthermore, we may restrict attention to the connected subgroup $H_0$. In fact, it is sufficient to establish the following: 
For each open $H_0$-orbit $\mathcal{O}'$, pick $\xi_0' \in \mathcal{O}'$ arbitrary. Assume that  there exists a map 
$\sigma: \mathcal{O}' \to H_0$ satisfying, for all $\xi \in \mathcal{O}'$,
\[
 \sigma(\xi)^T \xi_0' = \xi~,
\]
(i.e., $\sigma$ is a cross-section with respect to the dual action),
as well as
\[
 A(\xi)^e \| \sigma(\xi)^{\pm 1} \| \le C~,
\] with constants $e,C$ independent of  $\xi$. 

We then get for all $h \in H_0$ that
\[
 h \sigma(h^T \xi_0')^{-1} \in H_{\xi_0}~,
\]
and by assumption, the stabilizer is compact. Hence there exists a constant $K$ such that 
\[
 K := \sup_{h \in H} \left\| \left( h \sigma(h^T \xi_0')^{-1}\right)^{\pm 1} \right\| < \infty.
\]
But this yields, for all $h \in H_0$, that 
\begin{eqnarray*}
 A(h^T \xi_0')^{e} \| h^{\pm 1}\| & =  & A(\sigma(h^T \xi_0'))^e \| h^{\pm 1}\| \\
 & \le &  A(\sigma(h^T \xi_0')^e \| \sigma(h^T \xi_0')\| \left\| \left( h \sigma(h^T \xi_0')^{-1}\right)^{\pm 1} \right\| \\
 & \le & CK~,
\end{eqnarray*}
by assumption on $\sigma$. 

Now if $H$ is an irreducibly admissible matrix group containing $H_0$, then $H/H_0$ is finite, and the single open $H$-orbit is just the union of open $H_0$-orbits $\mathcal{O}'$. I.e., each 
$h \in H$ can be written as $h=g_i h_0$ with $g_i \in \{ g_1,\ldots,g_k\}$, the latter being a fixed system of representatives mod $H_0$, and $h_0 \in H_0$. Using a bound on the norms of the $g_i^{\pm 1}$, it is now easy to extend the desired estimate from $H_0$ to $H$. 

A similar reasoning allows to verify the condition of Proposition \ref{prop:van_moments}(a) by focusing on the connected component, as well; just replace all occurrences of $\| \cdot \|$ in the above argument by $|{\rm det}( \cdot)|$. 

In view of these observations, it remains to go through the list in \ref{thm:main}(a) and verify the conditions of Proposition \ref{prop:van_moments}. 
In the following, we will argue via \ref{prop:van_moments} (a) for case (2.l), and via \ref{prop:van_moments}(b) for the remaining cases, and obtain \ref{thm:main}(b),(c)  (recall Remark \ref{rem:atom_adm}).
 
\subsubsection{Known results: Cases (1) and (3)}

The existence of compactly supported atoms for abelian admissible groups acting in arbitrary dimensions was shown in \cite{Fu_RT}. This takes care of the cases (1.a) through (1.e). For the similitude group in arbitrary dimensions, this fact was established in \cite[Theorem 3.2]{Fu_atom}.

\subsubsection{The cases (2.a) through (2.i)}

In all these cases, one first verifies easily that 
\[
 \mathcal{O}_{\pm}  = H_0^T \left[ \begin{array}{c} \pm 1 \\ 0 \\ 0 \end{array} \right]  = \mathbb{R}^{\pm} \times \mathbb{R}^2
\] are the only open orbits, and we obtain the envelope function
\[
 A(\xi) = \min \left( \frac{|\xi_1|}{1+|[\xi_2,\xi_3]|}, \frac{1}{1+|\xi|} \right)~. 
\]

Furthermore, the dual action on the open orbits is free for all cases involved, and thus the cross-section $\sigma : \mathcal{O}_{\pm} \to H_0$ associated to $\xi_0 = [\pm 1, 0,0]^T$ is uniquely determined. 
Since the entries of $h^T [1,0,0]^T$ are just the entries of the first row of $h$, the computation of $\sigma(\xi)$ is equivalent to the task of computing the remaining entries of $h \in H_0$ from its first row entries. In each of the cases (2.d)-(2.i), one easily reads off the formulas for the group elements of $H_0$ that the entries of $\sigma(\xi)$ depend polynomially on $\xi_2, \xi_3$ and on arbitrary nonzero powers of $|\xi_1|$ and $\ln(|\xi_1|)$. Hence the estimates
\[
 A(\xi) |\xi_1|^{\pm 1} \le 1~, A(\xi) |\xi_{2}| \le 1~,A(\xi) |\xi_{3}| \le 1~
\] as well as $\lim_{x \to 0} x \ln(|x|) = 0$ yield boundedness of $A(\xi)^e \| \sigma(\xi)\| \le C$ for sufficiently large $e, C$. 

For the estimate of the norm of the inverse, we need the additional observation that ${\rm det}(\sigma(\xi)) = |\xi_1|^\alpha$ for a suitable real number $\alpha$, which follows from the fact that $\sigma(\xi)$ is upper triangular, with diagonal entries depending only on $\xi_1$. By Cramer's rule, the entries of $\sigma(\xi)^{-1}$ are polynomials in ${\rm det}(\sigma(\xi))^{-1}$ and the entries of $\sigma(\xi)$. Both can be controlled by suitable powers of $A(\xi)$, hence we get $A(\xi)^{e'} \| \sigma(\xi)^{-1} \| \le C'$ as well. 

The treatment of case (2.a) follows similar lines. Again, the action is free, and thus the cross-section is unique. Here, the moduli of the remaining entries of $\sigma(\xi)$ are less than or equal to $|\xi_1|$, and thus $A(\xi) \| \sigma(\xi) \| \le C$ follows. The estimate for the inverse again follows from this, Cramer's rule, and the dependence of ${\rm det}(\sigma(\xi))$ on $\xi_1$. 

For the remaining cases (2.b) and (2.c), note that the cross-section for (2.a) also works for these cases; to treat case (2.b), set $a=1$ in case (2.a). Hence we are done here as well.  

\subsubsection{The cases (2.j)-(2.k) and (2.m)}

For each of these groups, we obtain the four open orbits
\[
 \mathcal{O}_{\pm,\pm}  = H_0^T \left[ \begin{array}{c} \pm 1 \\ 0 \\ \pm 1 \end{array} \right]  = \mathbb{R}^\pm \times \mathbb{R} \times \mathbb{R}^\pm~,
\] and the associated envelope function
\[
 A(\xi) = \min \left( \frac{|\xi_1|}{1+|[\xi_2,\xi_3]|}, \frac{|\xi_3|}{1+|[\xi_1,\xi_2]|}, \frac{1}{1+|\xi|} \right)~. 
\]
Again the action is free, and thus the associated cross-section is unique. An arbitrary matrix $h \in H_0$ has the structure 
\[
 h = \left[ \begin{array}{ccc} h_{1,1} & h_{1,2} & 0 \\ 0 & h_{2,2} & 0 \\
             0 & 0 & h_{3,3} 
            \end{array} \right]
\]
and thus, with $\epsilon_1,\epsilon_2 \in \{ \pm 1 \}$, we get
\[
 h^T \left[ \begin{array}{c} \epsilon_1 \\ 0 \\ \epsilon_2 \end{array} \right]  = \left[ \begin{array}{c}\epsilon_1  h_{1,1} \\ \epsilon_1 h_{1,2} \\  \epsilon_2 h_{3,3} \end{array} \right]~.
\]
Thus the computation of $\sigma(\xi)$ boils down to computing $h_{2,2}$ from $h_{1,1},h_{1,2},h_{3,3}$, and in each of the cases under consideration, $h_{2,2}$ is seen to be a product of powers of $h_{1,1}$ and $h_{3,3}$. Now the estimates 
\[
 A(\xi) |\xi_{1}|^{\pm 1} \le 1~,  A(\xi) |\xi_{3}|^{\pm 1} \le 1~,A(\xi) |\xi_{2}| \le 1
\] yield $A(\xi)^e \| \sigma(\xi) \| \le C$ for suitably large $e, C$. The estimate for $\| \sigma(\xi)^{-1}\|$ again follows from this fact, Cramer's rule, and the fact that the determinant of $h$
is a product of powers of $h_{1,1}$ and $h_{3,3}$. 

\subsubsection{The case (2.n)}

We obtain the four open orbits
\[
 \mathcal{O}_{\pm,\pm}  = H_0^T \left[ \begin{array}{c} \pm 1 \\ \pm 1 \\ 0  \end{array} \right]  = \mathbb{R}^\pm \times \mathbb{R}^\pm \times \mathbb{R}~,
\] and the associated envelope function
\[
 A(\xi) = \min \left( \frac{|\xi_1|}{1+|[\xi_2,\xi_3]|}, \frac{|\xi_2|}{1+|[\xi_1,\xi_3]|}, \frac{1}{1+|\xi|} \right)~. 
\]
Again the action is free, and the associated cross-section is determined as
\[
 \sigma(\xi) = \left[ \begin{array}{ccc} |\xi_1| & 0 & \xi_3 -  \delta_1 \ln (|\xi_1/\xi_2|) - \delta_2 |\xi_2| \ln (|\xi_2|)  \\ 0 & |\xi_2| & \delta_1 \ln (|\xi_1/\xi_2|) + \delta_2 |\xi_2| \ln (|\xi_2|) \\ 0 & 0 & |\xi_2| \end{array} \right]~.
 \]
Again, all entries of $\sigma(\xi)$ are controlled by suitable positive powers of $|\xi_1|^{\pm 1}, |\xi_{2}|^{\pm 1}, |\xi_3|$, and $A(\xi)^e \| \sigma(\xi)\| \le C$ for suitably large $e,C$ follows. The estimate for the inverse now follows in the usual way from an appeal to Cramer's rule and the observation that $|{\rm det}(\sigma(\xi))| = |\xi_1|\, |\xi_2|^2$.  

\subsubsection{The exceptional case (2.l)}

We obtain the two open orbits
\[
 \mathcal{O}_{\pm,\pm}  = H_0^T \left[ \begin{array}{c} \pm 1 \\ 0 \\ 0  \end{array} \right]  = \mathbb{R}^\pm \times \mathbb{R}^2 ~,
\] and the associated envelope function
\[
 A(\xi) = \min \left( \frac{|\xi_1|}{1+|(\xi_2,\xi_3)|}, \frac{1}{1+|\xi|} \right)~. 
\]
This action is free, and the associated cross-section is determined as
\[
 \sigma(\xi) = \left[ \begin{array}{ccc} |\xi_1| & \xi_2 & \xi_3  \\ 0 & \exp(-\xi_3/\nu_1) |\xi_1|\exp(\nu_2 |\xi_1| \ln (|\xi_1|) / \nu_1)  & 0\\ 0 & 0 & |\xi_1| \end{array} \right] ~.
 \]
 Now one easily verifies that, for all $e \ge 0$, 
 \[ A(\sigma(1,0,\xi_3))^e \| \sigma(1,0,\xi_3) \| \ge C A(\sigma(1,0,\xi_3))^e \exp(-\xi_3/\nu_1) \to \infty~,
 \] as ${\rm sign}(\nu_1) \xi_3 \to -\infty$.  
 This shows that the required condition for vanishing moment criteria are in fact violated by this group. Hence the criteria for the existence of compactly supported atoms are not applicable. 
 
 We now show the existence of compactly supported admissible vectors. For this purpose, we need to compute the modular function of $H_0$. 
 Writing 
 \[
  h(a,b,t) = \left[ \begin{array}{ccc} b & t & -\nu_1 \ln a + \nu_2 b \ln b \\ 0 & ab & 0 \\ 0 & 0 & b \end{array} \right] \in H_0
 \] for $(a,b,t) \in \mathbb{R}^+ \times \mathbb{R}^+ \times \mathbb{R}$, we see that the subsets
 \[
  H_1 = \{ h(1,b,0): b \in \mathbb{R}^+ \}~, H_2 = \{ h(a,1,t): a \in \mathbb{R}^+, t \in \mathbb{R} \} 
 \] are two closed subgroups that commute elementwise. I.e., $H_0 \cong H_1 \times H_2$, and $H_1$ is isomorphic to $\mathbb{R}$. Furthermore, the computation
 \[
  h(a,1,t) h(a',1,t') = h(t'+a't,tt')
 \] shows that $H_2$ is (anti)-isomorphic to the well-known $ax+b$-group. Hence we get that the left Haar measure of $H_2$ is given by $dt \frac{da}{a}$, and the modular function of $H_2$ is $\Delta_{H_2}(h(a,1,t)) = |a|$.  In summary, this gives 
 \[
  \Delta_{H_0}(h(a,b,t)) = \Delta_{H_1}(h(1,b,0)) \Delta_{H_2}(h(a,1,t)) = |a|~,
 \]
 and thus 
 \[
  \Delta_{G}(h(a,b,t)) = \frac{\Delta_{H_0}(h(a,b,t))}{|{\rm det}(h(a,b,t))|} = \frac{|a|}{|a|\, |b|^3} =|b|^{-3}~.
 \]
But this implies for the function $\Phi$ entering the admissibility condition that 
\begin{eqnarray*}
 \Phi(\xi) & = & \Delta_G(\sigma(\xi_1,\xi_2,\xi_3)) = \Delta_G(h(\exp(-\xi_3/\nu_1) \exp(\nu_2 |\xi|_1 \ln (|\xi_1)| / \nu_1),\xi_1,\xi_2))  \\
 & = & |\xi_1|^{-3}~.
\end{eqnarray*}
In particular, $A(\xi)^3 \Phi(\xi) \le 1$, which proves the existence of compactly supported admissible vectors. 

\bibliographystyle{plain}

\bibliography{admissible_3d}

\end{document}